\documentclass [twoside,reqno,12pt] {amsart}

\usepackage{graphicx}

\usepackage{amsfonts}
\usepackage{amssymb}
\usepackage{a4}
\usepackage{color}

\newtheorem{thm}{Theorem}[section]
\newtheorem{cor}[thm]{Corollary}
\newtheorem{lem}[thm]{Lemma}
\newtheorem{prop}[thm]{Proposition}

\newtheorem{rem}[thm]{Remark}

\theoremstyle{definition}

\theoremstyle{remark}

\numberwithin{equation}{section}

\newcommand{\R}{\mathbb{R}}

\newcommand{\Z}{\mathbb{Z}}

\def\hat{\widehat}
\def\tilde{\widetilde}
\def \bfo {\begin {eqnarray*} }
\def \efo {\end {eqnarray*} }
\def \ba {\begin {eqnarray*} }
\def \ea {\end {eqnarray*} }
\def \beq {\begin {eqnarray}}
\def \eeq {\end {eqnarray}}

\def \det {\hbox{det}}

\def \p {\partial}

\def\hat{\widehat}
\def\tilde{\widetilde}
\def \bfo {\begin {eqnarray*} }
\def \efo {\end {eqnarray*} }
\def \ba {\begin {eqnarray*} }
\def \ea {\end {eqnarray*} }
\def \beq {\begin {eqnarray}}
\def \eeq {\end {eqnarray}}

\def \det {\hbox{det}}

\def \p {\partial}

\begin{document}

 \title[Inverse problems for differential forms]{
 Inverse problems for differential forms on  Riemannian manifolds with boundary}

\author[Krupchyk]{Katsiaryna Krupchyk}

\address
        {K. Krupchyk, Department of Mathematics and Statistics \\
         University of Helsinki\\
         P.O. Box 68 \\
         FI-00014   Helsinki\\
         Finland}

\email{katya.krupchyk@helsinki.fi}

\author[Lassas]{Matti Lassas}

\address
        {M. Lassas, Department of Mathematics and Statistics \\
         University of Helsinki\\
         P.O. Box 68 \\
         FI-00014   Helsinki\\
         Finland}

\email{matti.lassas@helsinki.fi}

\author[Uhlmann]{Gunther Uhlmann}

\address
       {G. Uhlmann, Department of Mathematics\\
       University of Washington\\
       Seattle, WA  98195-4350\\
       USA}
\email{gunther@math.washington.edu}

\maketitle

\begin{abstract} Consider a real-analytic orientable connected complete Riemannian manifold $M$ with boundary of dimension $n\ge 2$ and let $k$ be an integer $1\le k\le n$. In the case when $M$ is compact of dimension $n\ge 3$,  we show that the manifold and the metric on it can be reconstructed, up to an isometry, from the set of the Cauchy data for harmonic $k$-forms, given on an open subset of the boundary.  This extends a   result of \cite{LU01} when $k=0$. In the two-dimensional case,  the same conclusion is obtained when considering the set of the Cauchy data for harmonic $1$-forms.  
Under additional assumptions on the curvature of the manifold, we carry out the same program when $M$ is complete non-compact. In the case $n\ge 3$, this generalizes the results of  \cite{LTU03} when $k=0$.  In the two-dimensional case,  we are able to reconstruct the manifold from the set of the Cauchy data for harmonic $1$-forms. 
\end{abstract}

\section{Introduction and statement of results}
The purpose of this paper is to study the inverse problem of the determination of a complete Riemannian manifold $(M,g)$ with boundary from the Cauchy data of harmonic differential forms, given on an open subset of the  boundary.  We emphasize that the determination of a Riemannian manifold includes, in addition to the reconstruction of the metric, the  recovery of the topological and the differentiable structures of $M$. 

Motivated by the problem of electrical impedance tomography,  \cite{Nach88, SylUhl87},  
the issue of the reconstruction of a Riemannian manifold from the set of the Cauchy data of harmonic functions is a basic question in the field of inverse problems. Here the Cauchy data  can be represented  as the graph of the Dirichlet-to-Neumann map $\Lambda_g$, which is defined by solving  the Dirichlet problem
\[
\Delta_g u=0,\quad u|_{\p M}=f,
\]
with a given $f\in C^\infty(\p M)$, and setting $\Lambda_g (f)=\p_\nu u|_{\p M}$, where $\p_\nu$ is the exterior normal derivative.  In particular, the importance of  the Dirichlet-to-Neumann map $\Lambda_g$ is due to the fact that it encodes boundary measurements of voltage and current flux in electrical impedance tomography.

Since the works \cite{LTU03, LU01}, it is known that a complete  real-analytic connected Riemannian manifold of dimension $n\ge 3$ can be recovered, up to an isometry, from the knowledge of the map $\Lambda_g$, given on an open subset of the boundary.  It is natural to expect that the reconstruction of the manifold is also possible from the Dirichlet-to-Neumann map, associated to the Hodge Laplacian on differential forms. The first step in this direction has been made in \cite{JL05}, where it was shown that the full Taylor series of the metric tensor at the boundary can be recovered from the Dirichlet-to-Neumann map of the Hodge Laplacian on $k$-forms, $k=1,\dots, n$. In this work we complete this study by reconstructing the Riemannian manifold from such a Dirichlet-to-Neumann map. Our results are unconditional in the case when $M$ is compact, whereas in the complete non-compact case we require conditions on the curvature of the manifold, as well as a spectral assumption. 

In the case when the dimension of $M$ is equal to $2$, it was shown in   \cite{LU01} that only the conformal class of the compact real-analytic connected Riemannian manifold can be determined from the knowledge of $\Lambda_g$. This obstruction is due to the conformal invariance of the Laplacian on functions, i.e. 
\[
\Delta_{\sigma g}=\sigma^{-1}\Delta_g, \quad \sigma\in C^\infty(M), \quad \sigma>0.
\]
Furthermore, it was shown in \cite{LTU03} that there exist complete  two-dimensional Riemannian manifolds with boundary that are not conformally equivalent, but that have identical Dirichlet-to-Neumann maps $\Lambda_g$. 
The construction in \cite{LTU03} motivated the construction of invisibility cloaks in electrostatics \cite{GLU03}. See the review \cite{GKLU09}. We digress to review the construction in \cite{LTU03}.

Let $(M,g)$ be a compact 2-dimensional manifold with
non-empty boundary, let
$x_0\in M$ and consider manifold
\ba
\tilde M=M\setminus \{x_0\}
\ea
with the metric
\ba
\tilde g_{ij}(x)=\frac 1{d_M(x,x_0)^2}g_{ij}(x),
\ea
where $d_M(x,x_0)$ is the distance between $x$ and $x_0$ on $(M,g)$.
Then $(\tilde M,\tilde g)$ is a complete, non-compact 2-dimensional
Riemannian
manifold with the boundary $\p \tilde M=\p M$.
On the manifolds $M$ and $\tilde M$ we consider the boundary value
problems
\ba
\left\{\begin{array}{l}
\Delta_g u=0\quad \hbox{in $M$,}\\
u=f\quad \hbox{on $\p M$,}\end{array}\right. \quad\hbox{and}\quad
\left\{\begin{array}{l}
\Delta_{\tilde g} \tilde u=0\quad \hbox{in $\tilde M$,}\\
\tilde u=f\quad \hbox{on $\p \tilde M$,}\\
\tilde u\in L^\infty(\tilde M).\end{array}\right.
\ea
These boundary value problems are uniquely solvable and
define the Dirichlet-to-Neumann maps
\ba
\Lambda_{M,g}f=\p_\nu u|_{\p M},\quad
\Lambda_{\tilde M,\tilde g}f=\p_\nu \tilde u|_{\p \tilde M}.
\ea
As mentioned earlier in the two dimensional case
functions which are harmonic with respect to the metric $g$  stay harmonic
with respect to any
metric which is conformal to $g$,
one can see that
$\Lambda_{M,g}=\Lambda_{\tilde M,\tilde g}$.
This can be seen using e.g.\ Brownian motion or
capacity arguments. Thus, the boundary measurements for $(M,g)$ and $(\tilde M,\tilde g)$
coincide.
This gives a counter example for the inverse
electrostatic problem on Riemannian surfaces - even the topology
of possibly non-compact Riemannian surfaces can not be determined
using boundary measurements (see Figure 1).
\begin{figure}[htbp]
\begin{center}
\includegraphics[width=8cm]{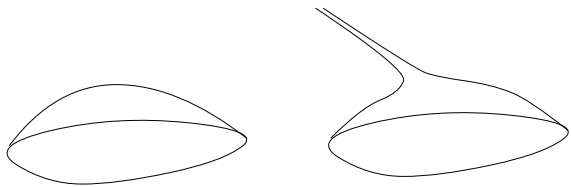} 
\caption{Blowing up a metric at a point, after \cite{LTU03}. The electrostatic
boundary measurements on the boundary of the  surfaces,
one compact and the other noncompact but complete, coincide.}
\end{center}
\end{figure}

In this paper, we exploit  the observation that the conformal invariance of the Laplacian can be broken by passing to forms of degree $1$.  Working on the level of $1$-forms, in the compact real-analytic case,  we reconstruct the manifold up to an isometry from the corresponding Dirichlet-to-Neumann map. In the complete non-compact case, we do the same under the assumption that the Gaussian curvature of the manifold should be bounded from below, as well as a spectral assumption.

The plan of the paper is as follows. The remainder of this section is devoted to the description of the precise assumptions, used throughout the paper, as well as to the statements of the main results. In Section 
\ref{sec_ex}, we construct an example of two complete  two-dimensional manifolds in the same conformal class, which can be told apart using the Cauchy data of harmonic $1$-forms.  The proofs of our results start in Section \ref{sec_compact}, where, following \cite{LTU03} closely, we address 
the case when $M$ is compact of dimension $n\ge 3$. The two-dimensional compact case is analyzed in Section 
\ref{sec_2dim}. Here an important role is played by the result that we can recover the Taylor series of the metric tensor at the boundary from the Dirichlet-to-Neumann map, associated to the Hodge Laplacian on $1$-forms.  The final Section 
\ref{sec_complete} establishes the results in the complete case.

\subsection{Notation}
Let $(M,g)$ be a real-analytic orientable connected complete  Riemannian  manifold of dimension $n\ge  2$ with a compact real-analytic boundary $\p M$ and let the metric $g$ be real-analytic up to the boundary.  Let $T^*M$ be the cotangent bundle on $M$ and let $\Lambda^k T^*M$, $k=0,1,\dots,n$,  be the bundles of the exterior differential $k$-forms. 
Denote by $C^\infty(M,\Lambda^k T^*M)$ the space of smooth  exterior differential forms of degree $k$. Here the smoothness is understood up to the boundary of $M$. 
The metric tensor $g$ induces the volume form $\mu=\mu_g\in C^\infty(M,\Lambda^n T^*M)$ and the Hodge  star isomorphism
\[
*:C^\infty(M,\Lambda^k T^*M)\to C^\infty(M,\Lambda^{n-k} T^*M), \quad \omega\wedge * \eta=g(\omega,\eta)\mu.
\]
Here in local coordinates, $\mu=\sqrt{\det (g_{ij})}dx^1\wedge \dots\wedge dx^n$, provided that $dx^1,\dots,dx^n$ is a positive basis of $T^*_x M$. 

Let $d: C^\infty(M,\Lambda^k T^*M)\to C^\infty(M,\Lambda^{k+1} T^*M)$ be the exterior  differential. Then 
 the codifferential operator is defined by
\[
\delta:C^\infty(M,\Lambda^k T^*M)\to C^\infty(M,\Lambda^{k-1} T^*M), \quad
\delta\omega=(-1)^{nk+n+1}*d*\omega,
\]
and the Hodge-Laplace operator is given by
\[
\Delta^{(k)}_{g}=\Delta^{(k)}:C^\infty(M,\Lambda^k T^*M)\to C^\infty(M,\Lambda^k T^*M), \quad \Delta^{(k)}\omega=(d\delta+\delta d)\omega.
\]
To study differential forms on the boundary of $M$, we consider the  inclusion map $i:\p M \to M$ and its pull-back,
\[
i^*:C^\infty (M,\Lambda^k T^*M)\to C^\infty (\p M,\Lambda^k T^*M).
\]
  Then we define the tangential trace of a $k$-form as
\[
\mathbf{t}:C^\infty (M,\Lambda^k T^*M)\to C^\infty (\p M,\Lambda^k T^*M), \quad \mathbf{t}\omega=i^*\omega, \quad  k=0,1,\dots,n-1,
\] 
and the normal trace as 
\[
\mathbf{n}:C^\infty (M,\Lambda^k T^*M)\to C^\infty (\p M,\Lambda^{n-k} T^*M), \quad \mathbf{n}\omega=i^*(*\omega), \quad k=1,2,\dots,n.
\]
The $L^2$-inner product of differential $k$-forms is given by
\[
(\omega,\eta)_{L^2}=\int_M
\omega\wedge * \bar \eta,\quad\omega,\eta\in C^\infty_0 (M,\Lambda^k T^*M),
\]
where $\bar \eta$ is the complex conjugate of $\eta$. The space $L^2(M,\Lambda^kT^*M)$ is defined as the completion of the space $C^\infty_0(M,\Lambda^kT^*M)$ of compactly supported differential $k$-forms on $M$ in the corresponding $L^2$-norm.  
When $\omega\in C^\infty_0(M,\Lambda^k T^* M)$, let
\[
\|\omega\|_s=(\sum_{j=0}^s\|\nabla^{(j)}\omega\|_{L^2}^2)^{1/2},
\]
where $\nabla^{(j)}$ is the $j$-th order covariant derivative of a $k$-form. 
The completion of the space $C^\infty_0(M,\Lambda^kT^*M)$ with respect to this norm is the standard Sobolev space of sections of $ \Lambda^kT^*M$, denoted by  $H^s(M, \Lambda^kT^*M)$  .

The space $C^\infty_0(M,\Lambda^k T^*M)$ is dense in $H^s(M,\Lambda^k T^*M)$ for $s\ge 0$, and the tangential and normal traces extend to continuous surjections
\begin{align*}
&\mathbf{t}: H^s (M,\Lambda^k T^*M)\to H^{s-1/2} (\p M,\Lambda^k T^*M),\\
&\mathbf{n}:H^s (M,\Lambda^k T^*M)\to H^{s-1/2} (\p M,\Lambda^{n-k} T^*M),
\end{align*}
as soon as $s>1/2$. The operators $d$ and $\delta$ extend to continuous mappings
\begin{align*}
d: H^s(M,\Lambda^k T^*M)\to H^{s-1}(M,\Lambda^{k+1} T^*M),\\
\delta:H^s(M,\Lambda^k T^*M)\to H^{s-1}(M,\Lambda^{k-1} T^*M),
\end{align*}
for $s\ge 1$. 

We
set
\[
\langle\mathbf{t}\omega,\mathbf{n}\eta\rangle=\int_{\partial
M}\mathbf{t}\omega\wedge\mathbf{n}\bar{\eta}, \quad\omega\in C^\infty (M,\Lambda^k T^*M),\ \eta\in C^\infty (M,\Lambda^{k+1} T^*M).
\]
From \cite[Proposition 2.1.2]{Sch_book} we recall the Stokes' formula
\begin{equation}
\label{eq_stokes_formulae}
(d\omega,\eta)_{L^2}-(\omega,\delta\eta)_{L^2}=
\langle\mathbf{t}\omega,\mathbf{n}\eta\rangle, \ \omega\in H^1 (M,\Lambda^k T^*M),\ \eta\in H^1 (M,\Lambda^{k+1} T^*M).
\end{equation}

The elements of the space
\[
\mathcal{H}^k(M)=\{\omega\in  H^1(M,\Lambda^kT^*M):d\omega=0,\delta\omega=0\}
\]
are called harmonic fields. Notice that the space $\mathcal{H}^k(M)$ is infinite dimensional for $1\le k\le n-1$, see \cite[Theorem 3.4.2]{Sch_book}. Two   subspaces are distinguished in $\mathcal{H}^k(M)$,
\begin{align*}
\mathcal{H}_D^k(M)&=\{\omega\in\mathcal{H}^k(M):\mathbf{t}\omega=0\}\quad\text{and}\\
\mathcal{H}_N^k(M)&=\{\omega\in\mathcal{H}^k(M):\mathbf{n}\omega=0\},
\end{align*}
which are called the Dirichlet and Neumann harmonic fields, respectively. According to \cite[Theorem 2.2.2]{Sch_book},  $\mathcal{H}_D^k(M)$ and $\mathcal{H}_N^k(M)$ are finite-dimensional, when $M$ is compact.

\subsection{The case of a compact manifold}  Here we assume that the manifold $M$ is compact with a non-empty boundary. Then the set of the Cauchy data of harmonic $k$-forms is given by
\begin{equation}
\label{eq_Cauchy}
\mathcal{C}_g^{(k)}=\{(\mathbf{t}\omega, \mathbf{n}\omega,\mathbf{n}d\omega,\mathbf{t}\delta \omega):\omega\in C^\infty (M,\Lambda^k T^*M), \Delta^{(k)}\omega=0\}.
\end{equation}
Here we notice that 
\[
\mathcal{C}_g^{(0)}=\{(\mathbf{t}\omega, \mathbf{n}d\omega):\omega\in C^\infty (M,\Lambda^0 T^*M), \Delta^{(0)}\omega=0\},
\]
and
\[
\mathcal{C}_g^{(n)}=\{( \mathbf{n}\omega,\mathbf{t}\delta \omega):\omega\in C^\infty (M,\Lambda^n T^*M), \Delta^{(n)}\omega=0\}.
\]
According to \cite{DS52,JL05}, the natural 
Dirichlet data for a harmonic $k$-form $\omega$ is 
\[
(\mathbf{t}\omega, \mathbf{n}\omega),
\]
and the natural Neumann data is 
\[
(\mathbf{n}d\omega,\mathbf{t}\delta \omega).
\]
Specifically, it is known \cite[Theorem 1]{DS52} that the problem 
\begin{equation}
\label{eq_dirichlet_problem}
\begin{aligned}
&\Delta^{(k)} \omega=0\quad \text{in}\quad  M, \\
&\mathbf{t}\omega=f_{1}\quad \text{on}\quad \p M,\\
&\mathbf{n}\omega=f_{2}\quad \text{on}\quad \p M.
\end{aligned}
\end{equation}
has a unique solution $\omega\in C^\infty (M,\Lambda^k T^*M)$ if $f_{1}\in C^\infty (\p M,\Lambda^k T^*M)$ and $f_2\in C^\infty (\p M,\Lambda^{n-k} T^*M)$. 
Thus, the natural Dirichlet-to-Neumann map for the $k$-form Hodge Laplacian is given by
\begin{equation}
\label{eq_DN_1}
\begin{aligned}
\Lambda^{(k)}_g:& \ C^\infty (\p M,\Lambda^k T^*M)\times C^\infty (\p M,\Lambda^{n-k} T^*M) \\
&\to C^\infty (\p M,\Lambda^{n-k-1} T^*M) \times C^\infty (\p M,\Lambda^{k-1} T^*M),\\
&\Lambda^{(k)}_g(f_1,f_2)=(\mathbf{n}d\omega,\mathbf{t}\delta \omega),
\end{aligned}
\end{equation}
where $\omega$ is the solution to \eqref{eq_dirichlet_problem}. 
Notice that the set of the Cauchy data \eqref{eq_Cauchy} for the harmonic $k$-forms is equal  to the graph of the Dirichlet-to-Neumann map \eqref{eq_DN_1}.

Another approach to the Dirichlet and Neumann boundary conditions for harmonic differential forms is given in  
\cite[Ch 5, Section 9]{Tay_bookI},
\begin{align*}
&\text{the Dirichlet boundary conditions}: \quad (\mathbf{t}\omega,\mathbf{t}\delta\omega),\\
&\text{the Neumann boundary conditions}: \quad (\mathbf{n}\omega,\mathbf{n}d\omega).
\end{align*}
It follows from \cite[Lemma 3.4.7]{Sch_book} that the problem 
\begin{equation}
\label{eq_dirichlet_2}
\begin{aligned}
&\Delta^{(k)} \omega=0\quad \text{in}\quad  M, \\
&\mathbf{t}\omega=g_1\quad \text{on}\quad \p M,\\
&\mathbf{t}\delta\omega=g_2\quad \text{on}\quad \p M,
\end{aligned}
\end{equation}
is solvable in $C^\infty (M,\Lambda^k T^*M)$, if and only if 
\[
\langle g_2,\mathbf{n}\lambda \rangle=0, \quad \forall \lambda\in\mathcal{H}^k_D(M).
\]
The solution of \eqref{eq_dirichlet_2} is unique up to an  arbitrary Dirichlet field $ \lambda\in\mathcal{H}^k_D(M)$. In this case, in order  to define the Dirichlet-to-Neumann map one has to specify the solution of \eqref{eq_dirichlet_2}, for instance by requiring that 
\begin{equation}
\label{eq_comp}
(\omega,\lambda)_{L^2}=0, \quad \forall \lambda\in\mathcal{H}^k_D(M).
\end{equation}
Notice that \eqref{eq_comp} gives a finite number of linear conditions. 
Let
\begin{align*}
W=\{(g_1,g_2)\in C^\infty (\p M,\Lambda^k T^*M)\times C^\infty (\p M,\Lambda^{k-1} T^*M):\langle g_2,\mathbf{n}\lambda \rangle=0,\\
 \forall \lambda\in\mathcal{H}^k_D(M)\}.
\end{align*}
Then the Dirichlet-to-Neumann map is defined by
\begin{align*}
&\tilde \Lambda^{(k)}_g:W\to C^\infty (\p M,\Lambda^{n-k} T^*M) \times C^\infty (\p M,\Lambda^{n-k-1} T^*M),\\
&\tilde \Lambda^{(k)}_g(g_1,g_2)=(\mathbf{n}\omega,\mathbf{n}d\omega),
\end{align*}
where $\omega$ is the solution of \eqref{eq_dirichlet_2} satisfying \eqref{eq_comp}.

\begin{prop}
\label{prop_cauchy}
The knowledge of  $\tilde \Lambda^{(k)}_g$ determines the set of the Cauchy data \eqref{eq_Cauchy} modulo a finite dimensional space. 
\end{prop}

\begin{proof}
Let 
\begin{align*}
U=\{&(\mathbf{t}\omega,\mathbf{n}\omega,\mathbf{n}d\omega,\mathbf{t}\delta\omega):\\
&\omega\in C^\infty (M,\Lambda^k T^*M),\Delta^{(k)}\omega=0,
\langle \mathbf{t}\delta\omega,\mathbf{n}\lambda\rangle=0,\forall \lambda\in \mathcal{H}^k_D(M)
\}.
\end{align*}
Then we can write
\[
\mathcal{C}^{(k)}_g=U\oplus\tilde U,
\]
where $\tilde U$ is finite-dimensional.  Furthermore,
\[
\mathcal{C}^{(k)}_g=\{(\mathbf{t}\omega,\mathbf{n}\omega,\mathbf{n}d\omega,\mathbf{t}\delta\omega)\in U:
(\omega,\lambda)_{L^2}=0, \forall \lambda\in\mathcal{H}^k_D(M)
\}\oplus\tilde U\oplus \hat U,
\]
$\hat U$ is finite-dimensional. The claim follows.
\end{proof}

Proposition \ref{prop_cauchy} implies that the knowledge of  $\tilde \Lambda^{(k)}_g$ determines   $\Lambda^{(k)}_g$ modulo a smoothing operator. 
Hence, it follows from \cite{JL05} that in case $n\ge 3$,
 the  knowledge of  $\tilde \Lambda^{(k)}_g$ allows us to recover the Taylor series at the boundary of the metric $g$ in the boundary normal coordinates.

Finally yet another definition of the Dirichlet-to-Neumann operator is given in 
\cite{BelSch1008}. For any $0\le k\le n-1$, the Dirichlet-to-Neumann operator 
\[
\hat \Lambda_g^{(k)}:C^\infty (\p M,\Lambda^k T^*M)\to C^\infty (\p M,\Lambda^{n-k-1} T^*M) 
\]
is defined as follows. Given $g_1\in C^\infty (\p M,\Lambda^k T^*M)$, the boundary problem
\begin{equation}
\label{eq_dirichlet_3}
\begin{aligned}
&\Delta^{(k)} \omega=0\quad \text{in}\quad  M, \\
&\mathbf{t}\omega=g_1\quad \text{on}\quad \p M,\\
&\mathbf{t}\delta\omega=0\quad \text{on}\quad \p M,
\end{aligned}
\end{equation}
is solvable in $C^\infty (M,\Lambda^k T^*M)$.  
The solution of \eqref{eq_dirichlet_3} is unique up to an arbitrary Dirichlet field $ \lambda\in\mathcal{H}^k_D(M)$.
Defining 
\[
\hat \Lambda_g^{(k)} g_1=\mathbf{n}d\omega,
\]
it is clear that $\hat \Lambda_g^{(k)}$ is independent of the choice of a solution of \eqref{eq_dirichlet_3}.

In \cite{BelSch1008}, an explicit formula is obtained which expresses the Betti numbers of the
manifold $M$ Êin terms of  $\hat \Lambda^{(k)}_g$.

The following is our first main result.
\begin{thm}
\label{thm_main}
Let $M_1$ and $M_2$ be compact real-analytic orientable connected Riemannian manifolds with real-analytic boundaries $\p M_1$ and $\p M_2$, respectively. Moreover, let the metrics $g_1$ and $g_2$ be real-analytic up to the boundary
and $\dim M_1=\dim M_2=n\ge 3$. 
Let $\Gamma_1=\Gamma_2=\Gamma$ be a non-empty open subset of $\p M_1$ and $\p M_2$. 
If for some integer  $1\le k\le n$,
\begin{equation}
\label{eq_C_G}
\mathcal{C}^{(k)}_{g_1}|_\Gamma=\mathcal{C}^{(k)}_{g_2}|_\Gamma,
\end{equation}
then $M_1$ and $M_2$ are isometric. 
\end{thm}

Here the sets $\Gamma_1\subset \p M_1$ and $\Gamma_2\subset \p M_2$ are identified by a real-analytic diffeomorphism.

Notice that if the manifolds $M_1$ and $M_2$ are isometric then the bundles of the exterior differential $k$-forms $\Lambda^k T^* M_1$ and $\Lambda^k T^* M_2$ are naturally isometric for any $k=1,\dots,n$. It follows therefore  from Theorem \ref{thm_main} that  the Cauchy data \eqref{eq_C_G} determine also the exterior differential form bundle structure over the manifold.

In the two dimensional case,  it is well-known \cite{LeeUhl89}  that there is an obstruction to reconstructing the metric from the set of Cauchy data for harmonic $0$-forms, due to the conformal invariance of the corresponding  Laplacian. It turns out that the conformal invariance can be broken by passing to forms of higher degree, and the reconstruction becomes possible. Thus, we obtain the following  result.

\begin{thm}
\label{thm_main_2dim}
Let $M_1$ and $M_2$ be compact real-analytic orientable connected Riemannian manifolds of dimension two with real-analytic boundaries $\p M_1$ and $\p M_2$, respectively. Moreover, let the metrics $g_1$ and $g_2$ be real-analytic up to the boundary. 
Let $\Gamma_1=\Gamma_2=\Gamma$ be a non-empty open subset of $\p M_1$ and $\p M_2$.  If  
\[
\mathcal{C}^{(1)}_{g_1}|_\Gamma=\mathcal{C}^{(1)}_{g_2}|_\Gamma,
\]
then $M_1$ and $M_2$ are isometric. 
\end{thm}

\subsection{The case of a complete non-compact  manifold}
Let $(M,g)$ be a complete non-compact real-analytic orientable connected Riemannian manifold with compact real-analytic boundary $\p M$. Let  
the metric $g$ be real-analytic up to the boundary and $n=\dim(M)\ge 2$.  

Let $f_1\in C^\infty(\p M,\Lambda^k T^*M)$ and $f_2\in C^\infty(\p M,  \Lambda^{n-k} T^* M)$. Then consider the problem
\begin{equation}
\label{eq_dirichlet_problem_complete}
\begin{aligned}
&\Delta^{(k)} \omega=0\quad \text{in}\quad  M, \\
&\mathbf{t}\omega=f_{1}\quad \text{on}\quad \p M,\\
&\mathbf{n}\omega=f_{2}\quad \text{on}\quad \p M.
\end{aligned}
\end{equation}

Let us make the following hypotheses.

\textbf{Assumption (A1)}. Let $n\ge 3$.  

\begin{itemize}

\item If $k=1$, then the  
Ricci curvature tensor is bounded from below.

\item If  $k=2,\dots, n$,  then the Riemann curvature tensor is bounded on $M$.

\end{itemize}

\textbf{Assumption (A2)}. Let $n=2$.  
Then the Gaussian curvature of $M$ is bounded from below.

In the two-dimensional case, similarly to the compact situation, we shall only be concerned with $1$-forms.

To discuss the solvability of the problem \eqref{eq_dirichlet_problem_complete} in the case when the complete manifold $M$ is non-compact, we shall need the following result. 

\begin{prop}
\label{prop_Friedrich} Assume that the assumption \emph{(A1)} holds. 
Then the  operator $\Delta^{(k)}=\Delta^{(k)}_F$, $k=1,\dots,n$, equipped with the domain
\begin{equation}
\label{eq_domain_1}
\mathcal{D}(\Delta^{(k)}_F)=\mathcal{K}\cap\{\omega\in L^2(M,\Lambda^k T^*M):\Delta^{(k)}\omega\in L^2(M,\Lambda^k T^* M)\},
\end{equation}
where
\[
\mathcal{K}=\{\omega\in H^1(M,\Lambda^k T^* M):\mathbf{t}\omega=0,\ \mathbf{n}\omega=0\},
\]
is a non-negative self-adjoint operator on $L^2(M,\Lambda^k T^* M)$. Assuming that the assumption \emph{(A2)} holds, we get the same conclusion for the operator $\Delta^{(1)}$. 

\end{prop}

\begin{proof}

The quadratic form 
\[
\mathcal{D}:\mathcal{K}\times \mathcal{K}\to \R, \quad \mathcal{D}(\omega,\omega)=\|d\omega\|^2_{L^2}+\|\delta \omega\|^2_{L^2},
\]
is densely defined and non-negative in $L^2(M,\Lambda^k T^* M)$. 

Let us show that $\mathcal{D}$ is closed, i.e. $\mathcal{K}$  is complete with respect to the norm
\[
|\!|\!|\omega|\!|\!|_{\mathcal{D}}=\sqrt{\mathcal{D}(\omega,\omega)+\|\omega\|^2_{L^2}}.
\] 
To this end, we shall use the following result from \cite[Theorem 2.1.5]{Sch_book}: as the boundary $\p M$ is compact, for all $\omega\in H^1(M,\Lambda^k T^*M)$ and $\mathbf{t}\omega=0$, we have
\begin{equation}
\label{eq_Caffney_almost}
\|\omega\|^2_{H^1}=\|\omega\|^2_{L^2}+(\mathcal{R}\omega,\omega)_{L^2}+\mathcal{D}(\omega,\omega)+\int_{\p M}\mathcal{S}\omega\wedge *\bar{\omega},
\end{equation}
where $\mathcal{R}\in \textrm{End}(\Lambda^k T^*M)$ is the bundle endomorphism, determined by the Riemannian curvature tensor on $M$ and $\mathcal{S}\in\textrm{End}(\Lambda^k T^*M|_{\p M})$ is the bundle endomorphism, determined by the second fundamental form of $\p M$.  

The compactness of $\p M$ implies that  the bundle endomorphism $\mathcal{S}$ is bounded, and 
\[
\big|\int_{\p M}\mathcal{S}\omega\wedge *\omega\big|\le C_{\mathcal{S}}\|\omega\|^2_{L^2(\p M)},\quad  C_{\mathcal{S}}>0.
\]
Furthermore, the compactness of $\p M$
 yields that for any $\epsilon>0$, there exists $C_\epsilon>0$ such that
\[
\|\omega\|^2_{L^2(\p M)}\le \epsilon \|\omega\|^2_{H^1(M)}+C_\epsilon\|\omega\|^2_{L^2( M)},
\]
see  \cite[Corollary 2.1.6]{Sch_book}.

If the Riemann curvature tensor is bounded on $M$, then the bundle endomorphism   $\mathcal{R}$ is also bounded, and
\[
\big| (\mathcal{R}\omega,\omega)_{L^2} \big|\le C_{\mathcal{R}} \|\omega\|^2_{L^2( M)}, \quad C_{\mathcal{R}}>0. 
\]
In the case of $1$-forms, $(\mathcal{R}\omega,\omega)_{L^2}=-\int_M\textrm{Ric}(\omega,\omega)\mu$, where $\textrm{Ric}$ is the Ricci tensor,
viewed as the bilinear form 
\[
\textrm{Ric}:\Lambda^1 T^*M\times \Lambda^1 T^*M\to \R.
\]
When $n=2$, the Ricci curvature tensor is  given by
\[
\textrm{Ric}^{ij}=Kg^{ij},
\]
where $K$ is the Gaussian curvature of $M$.

Hence, the assumptions of the proposition and \eqref{eq_Caffney_almost} imply that
\begin{equation}
\label{eq_Gaffney_2}
\|\omega\|_{H^1}\le C |\!|\!|\omega|\!|\!|_{\mathcal{D}},\quad C>0. 
\end{equation}
Thus, the fact that the quadratic form $D$ is closed follows from \eqref{eq_Gaffney_2} and the fact that 
the tangential $\mathbf{t}$ and normal $\mathbf{n}$ traces
\begin{align*}
&\mathbf{t}: H^1(M,\Lambda^kT^*M)\to H^{1/2}(\p M,\Lambda^k T^*  M), \quad k=0,1,\dots,n-1,\\
&\mathbf{n}: H^1(M,\Lambda^k T^*M)\to H^{1/2}(\p M,\Lambda^{n-k} T^* M), \quad k=1,2,\dots,n,
\end{align*}  
are continuous, see \cite[Theorem 1.3.7]{Sch_book}.  
The operator $\Delta^{(k)}_F$ is the non-negative self-adjoint operator, associated to the closed form $\mathcal{D}$, and the form domain of $\Delta^{(k)}_F$ is the space $\mathcal{K}$. 

\end{proof}

\begin{rem}
\label{rem_zero_forms}
The Hodge Laplace operator $\Delta^{(0)}$, acting on $0$-forms and equipped with the domain, given by \eqref{eq_domain_1}, is 
a non-negative self-adjoint operator on $L^2(M)$. This follows from the proof of Proposition \emph{\ref{prop_Friedrich}}, since in this case 
 the bundle endomorphism   $\mathcal{R}$ in \eqref{eq_Caffney_almost} vanishes, see \cite[Theorem 2.1.5]{Sch_book}.

\end{rem}

In what follows, assume that the assumption (A1) holds when $n\ge 3$, and the assumption (A2) holds when $n=2$. Denote by $\textrm{spec}_d(\Delta^{(k)}_F)$, $\textrm{spec}_{ess}(\Delta^{(k)}_F)$ and $\textrm{spec}(\Delta^{(k)}_F)$, the discrete spectrum, the  essential spectrum,  and the spectrum of  $\Delta^{(k)}_F$, respectively. 

\begin{rem}
\label{rem_zero_spec}
Let us show that $0\notin \textrm{spec}_d(\Delta^{(k)}_F)$. 
Take $\omega\in \mathcal{D}(\Delta^{(k)}_F)$ with $\Delta^{(k)}\omega=0$. Thanks to the boundary conditions, it follows that $d\omega=\delta\omega=0$. By an application of  \cite[Theorem 3.4.4]{Sch_book} , we get  $\omega=0$. 

\end{rem}

 In what follows, assume that $0\notin \textrm{spec}_{ess}(\Delta^{(k)}_F)$.  
To solve the problem \eqref{eq_dirichlet_problem_complete}, let us first observe that there is at most one solution $\omega\in (L^2\cap C^\infty)(M,\Lambda^k T^*M)$. Indeed,  if $\omega \in (L^2\cap C^\infty)(M,\Lambda^k T^*M)$, and $\Delta^{(k)}\omega=0$, $\mathbf{t}\omega=0$, $\mathbf{n}\omega=0$,  then 
by \cite[Theorem 2.1.5]{Sch_book}, $\omega\in \mathcal{D}(\Delta^{(k)}_F)$. Remark \ref{rem_zero_spec} implies that $\omega=0$.

To construct such a solution, let $F\in C^\infty(M, \Lambda^k T^*M)$ be an arbitrary form vanishing outside a bounded set,  such that 
\[
\mathbf{t}F=f_{1}, \quad \mathbf{n}F=f_{2}. 
\]
Since $0\notin \textrm{spec}(\Delta^{(k)}_F)$, there exists a unique  $v\in \mathcal{D}(\Delta^{(k)}_F)$ such that
\[
\Delta^{(k)}_F v=-\Delta^{(k)} F.
\]
Observe that by elliptic regularity, $v$ is $C^\infty$-smooth up to the boundary of $M$. We let 
$\omega=v+F$ be the unique solution of \eqref{eq_dirichlet_problem_complete} in $ (L^2\cap C^\infty)(M,\Lambda^k T^*M)$.

Therefore, the Dirichlet-to-Neumann map for the $k$-form Laplacian is determined by
\begin{align*}
\Lambda^{(k)}_g &:C^\infty(\p M,\Lambda^k T^*M)\times C^\infty(\p M,\Lambda^{n-k} T^* M)\\ 
&\to 
C^\infty(\p M,\Lambda^{n-k-1} T^* M)\times C^\infty(\p M,\Lambda^{k-1} T^* M),\\
&\Lambda^{(k)}_g(f_1,f_2)=(\mathbf{n}d\omega,\mathbf{t}\delta \omega),
\end{align*}
where $\omega\in  (L^2\cap C^\infty)(M,\Lambda^k T^*M)$ is the solution to \eqref{eq_dirichlet_problem_complete}.  
Let $\Gamma$ be a non-empty open subset of $\p M$ and $f_1\in C^\infty(\p M,\Lambda^k T^*M)$, $f_2\in C^\infty(\p M,\Lambda^{n-k} T^* M)$ be supported on $\Gamma$. Then we define the  Dirichlet-to-Neumann map related to $\Gamma$ by
\[
\Lambda^{(k)}_{g, \Gamma}(f_1,f_2)=(\mathbf{n}d\omega|_{\Gamma},\mathbf{t}\delta \omega|_\Gamma). 
\]
The following is our main result in the complete non-compact case. 

\begin{thm}
\label{thm_main_complete}
Let $M_1$ and $M_2$ be complete non-compact  real-analytic orientable connected Riemannian  manifolds,  $\dim(M_1)=\dim(M_2)=n\ge 2$, with compact real-analytic boundaries $\p M_1$ and $\p M_2$, respectively. Let the metrics $g_1$ and $g_2$ be real-analytic up to the boundary, and let $\Gamma_1=\Gamma_2=\Gamma$ be a non-empty open subset of $\p M_1$ and $\p M_2$. 
\begin{itemize}
\item If  $n\ge 3$, assume that \emph{(A1)} holds.  Suppose also that $0\notin \textrm{spec}_{ess}(\Delta^{(k)}_{F, M_1})$ and $0\notin \textrm{spec}_{ess}(\Delta^{(k)}_{F, M_2})$, and that  for some integer  $1\le k\le n$,
\[
\Lambda^{(k)}_{g_1, \Gamma}=\Lambda^{(k)}_{g_2, \Gamma}.
\]

\item If $n=2$, assume that \emph{(A2)} holds. Suppose also that $0\notin \textrm{spec}_{ess}(\Delta^{(1)}_{F, M_1})$ and $0\notin \textrm{spec}_{ess}(\Delta^{(1)}_{F, M_2})$, and that
\[
\Lambda^{(1)}_{g_1, \Gamma}=\Lambda^{(1)}_{g_2, \Gamma}.
\]
\end{itemize}
Then $M_1$ and $M_2$ are isometric. 
\end{thm}

\section{Examples in dimension two}
\label{sec_ex}

The purpose of this section is to illustrate Theorem \ref{thm_main_complete} by exhibiting an explicit example of two complete non-compact two-dimensional Riemannian manifolds, which can be distinguished from the knowledge of the set of the Cauchy data for harmonic $1$-forms,  while the sets of the Cauchy data for harmonic $0$-forms, produced by both manifolds, are identical.   

Let 
$M=S^1_{\theta}\times [0,+\infty)_z$. Here $S^1=\R/2\pi \Z$.  We shall consider two complete metrics on $M$, belonging to the same conformal class, 
\[
g_1=e^{-2z}(d\theta)^2+(dz)^2,
\]
and 
\[
g_2=(1+f(z))g_1,\quad f(z)=ze^{-z}.
\]
The Gaussian curvature of the manifold $(M,g_1)$ is  $K_{g_1}=-1$. 
We notice that $(M,g_1)$ is isometric to the following surface of revolution in $\R^3$, equipped with the standard metric,
\begin{equation}
\label{eq_sur_rev}
(z,\theta)\mapsto (h(z),e^{-z}\cos\theta,e^{-z}\sin\theta),
\end{equation}
where 
\[
h(z)=\int_0^z\sqrt{1-e^{-2t}}dt. 
\]
Explicitly, the surface of revolution \eqref{eq_sur_rev} is the familiar pseudosphere, obtained by rotating the curve
\begin{equation}
\label{eq_tractrix}
x(y)=-\sqrt{1-y^2}-\log y+\log(1+\sqrt{1-y^2}), \quad 0<y<1,
\end{equation}
about the $x$-axis, see Figure 2. 
\begin{figure}[htbp]
\begin{center}
\centerline{
\includegraphics[width=.4\linewidth]{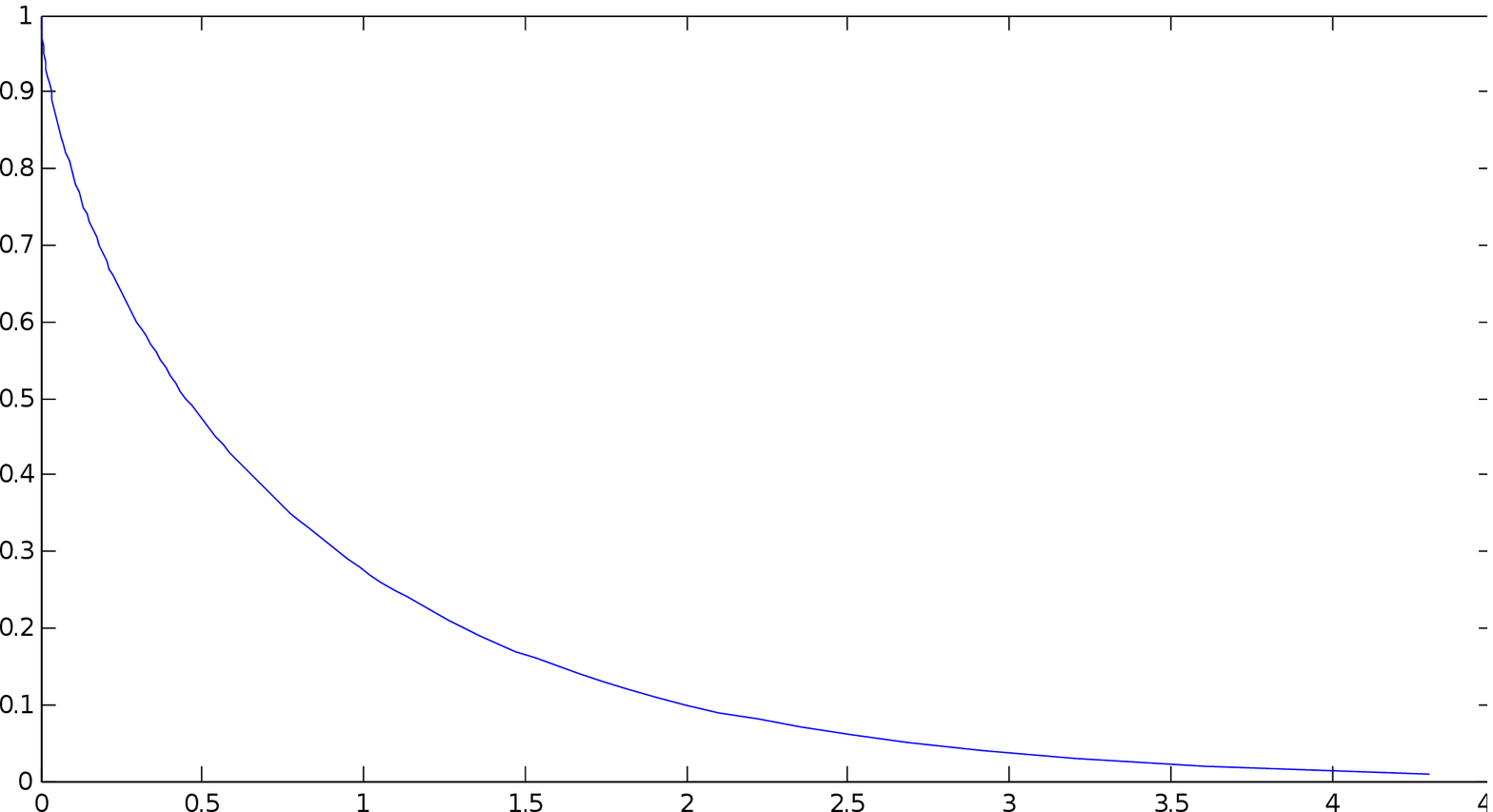}\hspace{10mm}
\includegraphics[width=.55\linewidth]{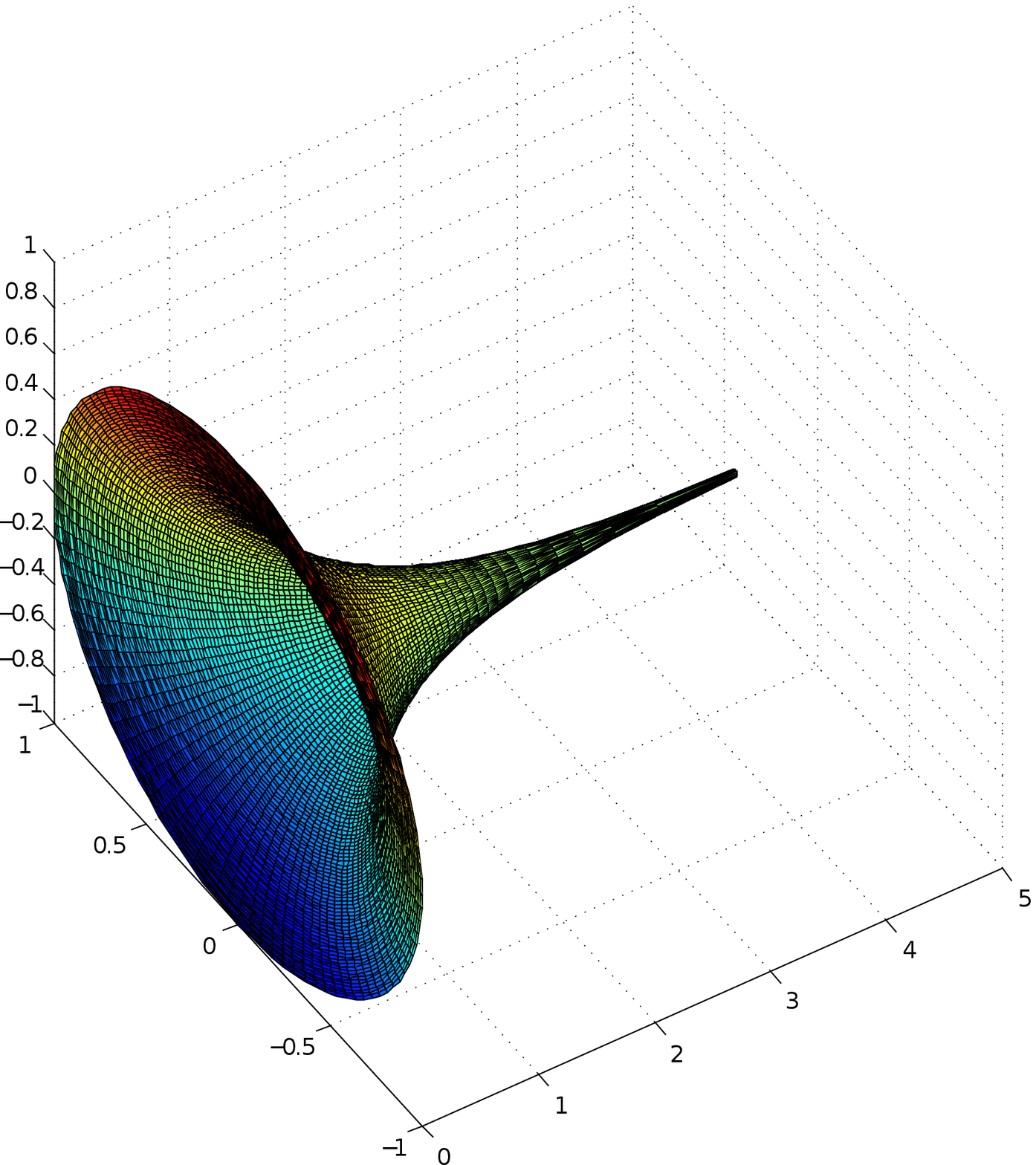}}
\end{center}
\caption{The surface of revolution on the right is obtained by rotating the tractrix \eqref{eq_tractrix}, on the left, about the $x$-axis.}
\end{figure}

 To compute the Gaussian curvature  $K_{g_2}$ of $(M,g_2)$, we use that
\[
\Delta^{(0)}_{g_1}(\frac{1}{2}\log (1+f(z)))+K_{g_1}=K_{g_2} (1+f(z)),
\]
see \cite{Bes_book}. Now
\begin{equation}
\label{eq_delta_0_ex}
\Delta^{(0)}_{g_1}=-e^{2 z}\p_\theta^2-\p_z^2+\p_z,
\end{equation}
and a direct computation gives that
\[
K_{g_2}=\frac{e^{-2z}(z+1-3z^2)+e^{-z}(3-6z)-2}{2(1+ze^{-z})^3}.
\]
In particular,  $|K_{g_2}|$ is uniformly bounded on $M$. Furthermore, since $K_{g_1}\ne K_{g_2}$,  the manifolds $(M,g_1)$ and $(M,g_2)$ are not isometric.

Let us now verify that zero is not in the spectrum of the Hodge Laplacian on $1$-forms $\Delta^{(1)}_{g_1}$, provided with the domain \eqref{eq_domain_1}. 
Writing $\omega(\theta,z)=\omega_1(\theta,z)d\theta+\omega_2(\theta,z)dz$ for a $1$-form $\omega$, a straightforward computation, using the formula \eqref{eq_lap_1} given below, 
 gives that 
\[
\Delta^{(1)}_{g_1}=\begin{pmatrix} -\p_z^2-e^{2z}\p_\theta^2-\p_z & 2\p_\theta\\
-2e^{2z}\p_\theta & -\p_z^2-e^{2z}\p_\theta^2+\p_z
\end{pmatrix}.
\]
This operator is considered in the Hilbert space 
\begin{equation}
\label{eq_L-2}
\mathcal{H}=L^2(S^1\times \R_+; e^{z}d\theta dz)\times L^2(S^1\times \R_+; e^{-z}d\theta dz),
\end{equation}
equipped with the inner product
\[
\langle \omega, \eta \rangle_{L^2}=\int_0^{2\pi} \int_0^{+\infty}(\omega_1\bar \eta_1e^z+\omega_2\bar \eta_2e^{-z})d\theta dz. 
\]
It follows from Proposition \ref{prop_Friedrich} that the form domain of the operator $\Delta^{(1)}_{g_1}$ is given by
\[
\mathcal{K}=\{\omega\in H^1(M,\Lambda^1 T^*M): \mathbf{t}w=0, \mathbf{n}w=0\}.
\]
Let us now obtain an explicit description of the space $\mathcal{K}$. Notice first that the boundary conditions take the form $\omega_1(\theta,0)=\omega_2(\theta,0)=0$. To make the condition that $\omega\in H^1(M,\Lambda^1 T^*M)$ explicit, we shall use \eqref{eq_Caffney_almost}. We get
\begin{equation}
\label{eq_cond_h1}
\begin{aligned}
&\int_0^{2\pi}\int_0^{+\infty}(|\omega_1|^2e^z+|\omega_2|^2e^{-z})d\theta dz< \infty,\\
& \int_0^{2\pi}\int_0^{+\infty} (e^{3z}|\p_\theta\omega_1|^2+e^z |\p_z\omega_1|^2)d\theta dz< \infty,\\
&\int_0^{2\pi}\int_0^{+\infty}(e^{z}|\p_\theta\omega_2|^2+e^{-z} |\p_z\omega_2|^2)d\theta dz< \infty.
\end{aligned}
\end{equation}
It is now convenient to perform a conjugation by the unitary operator  
\[
U=\begin{pmatrix} e^{z/2} & 0\\
0& e^{-z/2}
\end{pmatrix},
\]
in order to pass to the unweighted space $(L^2(S^1\times \R_+; d\theta dz))^2$. 
We have
\[
U\Delta^{(1)}_{g_1} U^{-1}
=\begin{pmatrix} -\p_z^2 +1/4- e^{2z} \p_\theta^2 & 2 e^z \p_\theta\\
-2 e^z \p_\theta & -\p_z^2+1/4-e^{2z} \p_\theta^2
\end{pmatrix},
\]
and the conditions \eqref{eq_cond_h1}, when expressed in terms of 
\[
u_1=e^{z/2}\omega_1\in L^2(S^1\times \R_+,d\theta dz),\quad u_2=e^{-z/2}\omega_2\in L^2(S^1\times\R_+, d\theta dz),
\]
become 
\begin{equation}
\label{eq_h1-2}
\int_0^{+\infty}\int_0^{2\pi}(e^{2z}|\p_\theta u_j|^2+ |\p_z u_j|^2)dzd\theta< \infty, \quad j=1,2. 
\end{equation}
Taking a Fourier decomposition with respect to the variable $\theta\in S^1$, we write
\begin{equation}
\label{eq_P_k}
U\Delta^{(1)}_{g_1} U^{-1}=\bigoplus_{k\in \Z}P_k,\quad P_k=\begin{pmatrix} -\p_z^2 +1/4+e^{2z} k^2 & 2ik e^z\\
-2ik e^z & -\p_z^2+1/4+e^{2z} k^2
\end{pmatrix}.
\end{equation}
It follows from \eqref{eq_h1-2} that the form domain of the operator $P_k$, $k\ne 0$, is given by
\begin{equation}
\label{eq_form_dom_p_k}
\{(u_1,u_2)\in (L^2(\R_+))^2: \int_0^{+\infty}(|\p_z u_j|^2+ e^{2z}|u_j|^2)dz<\infty, u_j(0)=0, j=1,2 \}.
\end{equation}
Let us check that this space is compactly embedded in $(L^2(\R_+))^2$. Indeed, let $(f_l)_{l=1}^\infty$ be such that
\[
\int_0^{+\infty}(|\p_z f_l|^2+ e^{2z}|f_l|^2)dz\le 1, \quad l=1,2,\dots. 
\]
Then
\begin{align*}
\|f_l-f_m\|^2_{L^2}&\le \int_0^R |f_l-f_m|^2dz+e^{-R}\int_R^{+\infty} e^z|f_l-f_m|^2dz\\
&\le \int_0^R |f_l-f_m|^2dz+4 e^{-R},\quad R>0. 
\end{align*}
Since the embedding $H^1((0,R))\to L^2((0,R))$ is compact for any $R>0$, the claim follows. 
Hence, the contribution to the continuous spectrum of  the operator $\Delta^{(1)}_{g_1}$ with the Dirichlet boundary conditions comes only from $k=0$ in \eqref{eq_P_k}.   
Here we have the operator $P_0$, whose spectrum is $[1/4,+\infty)$ and is purely absolutely continuous. 
We conclude that zero is not in the spectrum of the self-adjoint realization of $\Delta^{(1)}_{g_1}$, constructed in Proposition \ref{prop_Friedrich}.  See also \cite{DonXav84}.

Let us now also show that zero is also not in the spectrum of $\Delta^{(1)}_{g_2}$, provided with the domain \eqref{eq_domain_1}. To compute the coordinate expression for $\Delta^{(1)}_{g_2}$, we shall 
take advantage of the fact that the metric $g_2$ is a conformal multiple of $g_1$ and 
use the following identity,
\[
\Delta^{(1)}_{g_2}\omega=\tau \Delta_{g_1}^{(1)}\omega-i_{\nabla_{g_1}\tau}d\omega+d\tau\wedge \delta_{g_1}\omega,
\] 
where 
\[
\tau(z)=\frac{1}{1+f(z)}
\]
 and $i_{X}$ is the contraction  of a differential form with a vector field $X$, see \cite{Bes_book}.
We get
\[
\Delta^{(1)}_{g_2}=\tau(z)\Delta^{(1)}_{g_1}-\tau'(z)\begin{pmatrix}
\p_z & -\p_\theta\\
e^{2z}\p_\theta & \p_z-1
\end{pmatrix}.
\]
Since $*_{g_2}=*_{g_1}$ on $1$-forms, it follows that 
the operator $\Delta^{(1)}_{g_2}$ should be considered in the Hilbert space $\mathcal{H}$ in
\eqref{eq_L-2}. 
As $1/2<\tau(z)\le 1$, an explicit computation shows that the form domain of $\Delta^{(1)}_{g_2}$  is the same as the form domain of $\Delta^{(1)}_{g_1}$, i.e. 
\[
\{(\omega_1(\theta,z),\omega_2(\theta,z))\in\mathcal{H}: \eqref{eq_cond_h1}\ \textrm{holds},\ \omega_1(\theta,0)=\omega_2(\theta,0)=0
\}.
\]
Conjugating by the unitary operator $U$, as before, we get 
\begin{align*}
U\Delta^{(1)}_{g_2} U^{-1}
 &= \tau(z)\begin{pmatrix} -\p_z^2 +1/4-e^{2z} \p_\theta^2 & 2 e^z\p_\theta\\
-2 e^z\p_\theta & -\p_z^2+1/4-e^{2z}\p_\theta^2 
\end{pmatrix}\\
&-\tau'(z)
\begin{pmatrix}
\p_z-1/2 & - e^z\p_\theta\\
 e^z\p_\theta & \p_z-1/2 
\end{pmatrix}
\end{align*}
\small
\[
=\begin{pmatrix} -\p_z(\tau(z)\p_z) +\frac{\tau(z)}{4}+\frac{\tau'(z)}{2}-\tau(z)e^{2z}\p_\theta^2 & (2\tau(z)+\tau'(z))e^z\p_\theta\\
-(2\tau(z)+\tau'(z))e^z\p_\theta&  -\p_z(\tau(z)\p_z) +\frac{\tau(z)}{4}+\frac{\tau'(z)}{2}-\tau(z)e^{2z}\p_\theta^2
\end{pmatrix}.
 \]
\normalsize
Taking the Fourier decomposition with respect to the variable $\theta\in S^1$, we have
\[
U\Delta^{(1)}_{g_2} U^{-1} =\bigoplus_{k\in \Z} L_k, 
\]
where $L_k$ is defined by
\small
\[
\begin{pmatrix} -\p_z(\tau(z)\p_z) +\frac{\tau(z)}{4}+\frac{\tau'(z)}{2}+\tau(z)e^{2z}k^2 & (2\tau(z)+\tau'(z))e^zik\\
-(2\tau(z)+\tau'(z))e^zik&  -\p_z(\tau(z)\p_z) +\frac{\tau(z)}{4}+\frac{\tau'(z)}{2}+\tau(z)e^{2z}k^2
\end{pmatrix}.
 \]
\normalsize
The form domain of the operator $L_k$, $k\ne 0$, is given by \eqref{eq_form_dom_p_k}, and it is compactly embedded in $(L^2(\R_+))^2$. Hence, the spectrum of $L_k$, $k\ne 0$, is discrete. When $k=0$, we have
\[
L_0=\bigg(-\p_z(\tau(z)\p_z) +\frac{\tau(z)}{4}+\frac{\tau'(z)}{2} \bigg)\otimes I.
\]
Since the potential $\tau(z)/4+\tau'(z)/2\to 1/4$ as $z\to +\infty$, it follows from 
\cite[Chapter 8]{Dav_book} that the essential spectrum of $L_0$ is contained in $[1/4,+\infty)$. 
We conclude that $0$ is not in the spectrum of $\Delta^{(1)}_{g_2}$.

An application of Theorem \ref{thm_main_complete} shows therefore that the set of the Cauchy data for harmonic $1$-forms allows us to  distinguish between the Riemannian manifolds $(M,g_1)$ and $(M,g_2)$.

We shall conclude this example by pointing out that the manifolds $(M,g_1)$ and $(M,g_2)$ are undistinguishable on the level of the Cauchy data for harmonic $0$-forms.   Recall from Remark \ref{rem_zero_forms} that the Hodge Laplace operator $\Delta^{(0)}_{g_1}$ and $\Delta^{(0)}_{g_2}$, equipped with the domain, given by \eqref{eq_domain_1}, are non-negative self-adjoint operators on $L^2(M)$. 
Using \eqref{eq_delta_0_ex}, we see as above that zero is not in the spectrum of $\Delta^{(0)}_{g_1}$. 
As 
\begin{equation}
\label{eq_delta_0_new}
\Delta^{(0)}_{g_2}=(1+f(z))^{-1}\Delta^{(0)}_{g_1},
\end{equation}
the same conclusion holds for the spectrum of $\Delta^{(0)}_{g_2}$. It follows that given $v\in C^\infty(\p M)$, the boundary value problem
\begin{equation}
\label{eq_lap_ex1}
\begin{aligned}
\Delta^{(0)}_{g_j} w_j&=0,\quad \textrm{on}\ M,\\
\mathbf{t}w_j&=v,
\end{aligned}
 \end{equation}
 has a unique solution $w_j\in L^2(M,\mu_{g_j})$, $j=1,2$. Since $\mu_{g_2}=(1+f(z))\mu_{g_1}$, we have
 $L^2(M,\mu_{g_2})=L^2(M,\mu_{g_1})$ as spaces, with equivalent norms. In view of \eqref{eq_delta_0_new}, we get $w_1=w_2=w$. 
Moreover, $\mathbf{n}_{g_2}dw=\mathbf{n}_{g_1}dw$. It follows that the manifolds $(M,g_1)$ and $(M,g_2)$ produce the same set of the Cauchy data of  harmonic $0$-forms. 

\begin{rem}
In the paper \cite{LTU03}, the case of the Laplacian on $0$-forms on a complete non-compact manifold was considered. The construction of the solution to 
\eqref{eq_lap_ex1} in  \cite{LTU03} is based on the maximum principle and leads to a solution in $L^\infty(M)$, whereas in the present paper $L^2$-methods are employed.  The definition of the Dirichlet-to-Neumann map in   \cite{LTU03} is consequently different from the definition used in this paper. Nevertheless, it is clear that the sets of the Cauchy data of the manifolds $(M,g_1)$ and $(M,g_2)$, defined using the $L^\infty$-approach, are identical.

\end{rem}

\section{The case of a compact manifold. Proof of Theorem \ref{thm_main}}

\label{sec_compact}

\subsection{Reconstruction near the boundary}

 Near $\Gamma$ let us consider boundary normal coordinates given by a local coordinate chart $(x^1,\dots,x^n)$, where $x^n\ge 0$ is the distance to $\Gamma$ and $x'=(x^1,\dots,x^{n-1})$ is a local chart on  $\Gamma$.  In these coordinates the metric $g$ has the form, see e.g. \cite{LeeUhl89},
\[
g=\sum_{i,j=1}^{n-1}g_{ij}(x)dx^idx^j +(dx^n)^2.
\] 
It was proved in \cite{JL05} that 
for any integer $k$, $1\le k\le n$, the set of Cauchy data $\mathcal{C}^k_g|_{\Gamma}$ determines all the normal derivatives $\p_{x^n}^l g_{ij}(x',0)$, $l=0,1,2,\dots$, of the metric tensor $g$ at $\Gamma$. 

By considering the Taylor expansion of the metric $g$ near a given  point on $\Gamma$, we extend the metric to a boundary collar  $\Gamma\times (-r,0]$, for $r>0$ small enough, so that the extended metric remains  real-analytic.  This metric is also uniquely determined. 
We introduce the real-analytic  manifold $\tilde M$ obtained by attaching to $M$ a boundary collar $\Gamma\times (-r,0]$, equipped with this metric.

\subsection{Green's forms}

\label{sec_Green's_forms}

It is known   \cite[p.139]{DS52} that for any $k=1,\dots,n$, there is a Green's form $G(x,y)$, which is a double form of degree $k$,  satisfying 
\begin{align*}
&\Delta_{x}^{(k)}G(x,y)=\delta_{x,y}
\quad \text{in}\quad \tilde M,\\
&\mathbf{t}_xG(x,y)=0\quad \text{on}\quad \p\tilde M,\\
&\mathbf{n}_xG(x,y)=0\quad \text{on}\quad \p\tilde M,
\end{align*}
where $y\in \tilde M\setminus\p\tilde M$ and 
\[
\delta_{x,y}=\sum_{\underset{
j_1<\dots <j_k}{i_1<\dots <i_k}}\delta(y-x)dx^{i_1}\wedge \dots\wedge dx^{i_k}\cdot
dy^{j_1}\wedge \dots\wedge dy^{j_k}
\]
 is the delta double current, supported at $x=y$. 
 See \cite{Rham_book} for the notions of a double form and a double current. 
 
  Throughout the following discussion the degree $k$, $k=1,\dots,n$, will be kept fixed. 
We have, $G(x,y)=G(y,x)$. Furthermore, as $x\to y$, $G(x,y)$ has the following asymptotic behavior  (see  \cite[p.136]{DS52})
\begin{equation}
\label{eq_asymptotic}
G(x,y)\sim\frac{1}{(n-2)s_nd_{\tilde M}(x,y)^{n-2}} \Gamma_{(i_1,\dots,i_k),(j_1,\dots,j_k)}(y)dx^{i_1}\wedge \dots \wedge dx^{i_k}\cdot dy^{j_1}\wedge \dots \wedge dy^{j_k},
\end{equation}
where $d_{\tilde M}(x,y)$ is the geodesic distance from $x$ to $y$ in $\tilde M$, $s_n$ is the area of the unit $n$-sphere and 
\[ 
\Gamma_{(i_1,\dots,i_k),(j_1,\dots,j_k)}(y)=\begin{vmatrix}
 g_{i_1j_1}(y) & \dots & g_{i_kj_1}(y)\\
 \dots & \dots &\dots\\
 g_{i_1j_k}(y) & \dots & g_{i_kj_k}(y)
 \end{vmatrix}.
\]
Here $\Gamma_{(i_1,\dots,i_k),(j_1,\dots,j_k)}(y)$ is positive real-analytic function on $\tilde M$, as it is a Gram determinant. 
The asymptotic relation \eqref{eq_asymptotic} can be differentiated  with respect to $x$ any number of times. 

\begin{lem} Let $y\in \tilde M\setminus\p \tilde M$. Then the Green's form $G(x,y)$ is real-analytic for $x\in \tilde M\setminus(\{y\}\cup \p\tilde M)$. 
\end{lem}

This lemma follows from \cite[Theorem 4.1, p. 108]{EgoShu}, since $G(x,y)$ is a solution of  an  elliptic system of partial differential equations on $\tilde M\setminus(\{y\}\cup \p\tilde M)$ with real-analytic coefficients. 

By the symmetry of the Green's forms, we have the following 

\begin{cor} Let $x\in \tilde M\setminus\p \tilde M$. Then the Green's form $G(x,y)$ is real-analytic for $y\in \tilde M\setminus(\{x\}\cup \p\tilde M)$.
\end{cor}

We shall also need the following fact which follows from \cite[Theorem 1]{DS52}.

\begin{lem} 
\label{lem_smooth}
Let $y\in \tilde M\setminus\p \tilde M$. Then the Green's form $G(x,y)$ is $C^\infty$-smooth up to the boundary of $\p \tilde M$, away from $y$. 
\end{lem}

\subsection{Reconstruction of the Green's forms in the collar neighborhood of the boundary}

Recall that near $\Gamma$ we consider a coordinate chart $(x^1,\dots,x^n)$ where $x^n$ is the distance to $\Gamma$ and $x'=(x^1,\dots,x^{n-1})$ is a chart on  $\Gamma$. We write a $k$-form as
\[
\omega=\sum_{i_1<\dots <i_k}\omega_{i_1,\dots,i_k}dx^{i_1}\wedge\dots\wedge dx^{i_k}.
\]
The restriction of $\omega$ to $\Gamma$ is given by
\[
\omega|_{\Gamma}=\sum_{i_1<\dots <i_k}\omega_{i_1,\dots,i_k}|_{\Gamma}dx^{i_1}\wedge\dots\wedge dx^{i_k}.
\]
The  restriction of the normal derivative of $\omega$ to $\Gamma$ is defined by
\[
\p_{x^n} \omega|_{\Gamma}=\sum_{i_1<\dots <i_k}\p_{x^n}\omega_{i_1,\dots,i_k}|_{\Gamma}dx^{i_1}\wedge\dots\wedge dx^{i_k}.
\]
Observe, as explained in \cite{JL05},
 that given the induced metric on $\Gamma$ and the first normal derivative of the metric on $\Gamma$, the knowledge of the Dirichlet data $(\mathbf{t}\omega,\mathbf{n}\omega)$ on $\Gamma$ implies the knowledge of $\omega|_{\Gamma}$ and the knowledge of the Neumann data $(\mathbf{t}\delta\omega,\mathbf{n}d\omega)$ on $\Gamma$ determines  $\p_{x^n} \omega|_{\Gamma}$. 

For future references, let us also recall that the Hodge Laplacian on $k$-forms
is a second order elliptic partial differential system, with a scalar principal symbol, which 
 in local coordinates given by 
 \begin{equation}
\label{eq_laplace}
\Delta^{(k)}= -(g^{jk}(x)\p_{x^j}\p_{x^k})\otimes I + B_j(x)\p_{x^j}+C(x),
\end{equation}
where $(g^{jk})=(g_{jk})^{-1}$, $I$ is  $d\times d$ identity matrix, and $B_j(x)$, $C(x)$ are real-analytic functions with values in the set of $d\times d$-matrices, $d=\begin{pmatrix} n\\ k\end{pmatrix}$.

Denote $U= \tilde M_1\setminus M_1= \tilde M_2\setminus M_2$. Strictly speaking, here the manifolds $\tilde M_1\setminus M_1$  and $\tilde M_2\setminus M_2$ are identified by a real-analytic isometry. 

\begin{lem}
\label{lem_green_new}
The Green's forms $G_j(x,y)$ satisfy
\[
G_1(x,y)=G_2(x,y), \quad (x,y)\in U\times U\setminus\{x=y\}.
\]
\end{lem}

\begin{proof}
Let $y\in   U$ and consider the following Dirichlet problem on $M_2$,
\begin{align*}
&\Delta_{g_2}^{(k)} \omega(x)=0,\quad x\in  M_2, \\
&\mathbf{t}\omega(x)=\mathbf{t}_x G_1(x,y),\quad x\in  \Gamma,\\
&\mathbf{n}\omega(x)=\mathbf{n}_x G_1(x,y),\quad x\in \Gamma,\\
&\mathbf{t}\omega(x)=0,\quad x\in  \p M_2\setminus\Gamma,\\
&\mathbf{n}\omega(x)=0,\quad x\in  \p M_2\setminus\Gamma.
\end{align*}
Notice that in this problem the boundary data are continuous. By \cite[Theorem 3.4.10]{Sch_book} the problem above is uniquely solvable, and the hypothesis that $\mathcal{C}_{g_2}|_\Gamma=\mathcal{C}_{g_1}|_\Gamma$ implies that
\begin{align*}
&\mathbf{n}d\omega(x)=\mathbf{n}_xdG_1(x,y),\quad x\in \Gamma,\\
&\mathbf{t}\delta \omega(x)=\mathbf{t}_x\delta G_1(x,y), \quad x\in \Gamma.
\end{align*}

Let 
\[
\tilde \omega(x)=\begin{cases} \omega(x), & x\in  M_2, \\
G_1(x,y), & x\in  U. 
\end{cases}
\]
Now since the Cauchy data of $\omega(x)$ and $G_1(x,y)$ coincide on $\Gamma$,   it follows from \eqref{eq_laplace} that $\tilde \omega$ satisfies
\begin{align*}
&\Delta_{g_2}^{(k)}\tilde \omega(x)=\delta_{x,y},\quad x\in \tilde M_2,\\
& \mathbf{t}\tilde \omega(x)=0,\quad x\in \p\tilde M_2,\\
&\mathbf{n}\tilde \omega(x)=0,\quad x\in \p\tilde M_2.
\end{align*}
Thus, $\tilde\omega(x)-G_2(x,y)$ solves the homogeneous Dirichlet problem on $\tilde M_2$. Since the manifold $\tilde M_2$ is compact, the solution to the Dirichlet problem is unique, cf. \cite[Theorem 1]{DS52}, i.e. $\tilde\omega(x)=G_2(x,y)$ for $x\in \tilde M_2\setminus\{y\}$.

\end{proof}

\subsection{Embedding of the manifold into a Sobolev space} 

Following \cite{LTU03} we will prove Theorem \ref{thm_main} via certain embeddings of $\tilde M_j$ into a Sobolev space given by the Green's forms. 
Let $\tilde U\subset\subset U$ be a relatively compact subset of $U$. Then 
for some $s<1-n/2$, we define the maps
\begin{equation}
\label{eq_embedding}
\mathcal{G}_j:\tilde M_j\to H^s(\tilde U,\Lambda^k T^* \tilde U),\quad
 x\mapsto  G_j(x,\cdot).
\end{equation}

\begin{prop}
\label{prop_regularity_emb}
The maps $\mathcal{G}_j$ are $ C^1(\tilde M_j,H^s(\tilde U,\Lambda^kT^* \tilde U))$, for $s<1-n/2$. 
\end{prop}

\begin{proof}
Let $V$ be a relatively compact subset of $U$ such that
\[
\tilde U\subset\subset V\subset \subset U.
\]
Assume first that $x\in \tilde M_j\setminus V$. Then it is known \cite[Chapter 7]{Tay_bookII} that
\[  
G_j(x,y)\in C^\infty(\tilde M_j\setminus V\times \overline{\tilde U}, \Lambda^k T^*(\tilde M_j\setminus V\times \overline{\tilde U})).
\]
Hence,
\begin{equation}
\label{eq_G_1}
\mathcal{G}_j |_{\tilde M_j\setminus V}\in C^\infty(\tilde M_j\setminus V, H^l(\tilde U,\Lambda^k T^* \tilde U)), \quad\textrm{for any}\ l. 
\end{equation}
Let now $x\in V$. It is easily seen that the following map is $C^1$,
\[
V\ni x\mapsto \delta_{x,y}\in H^{s-2}(V,\Lambda^k T^* V),
\]
since we have assumed that $s<1-n/2$.
As $G_j(x,y)=G_j(x,y)$, we get
\[
\Delta_y^{(k)}G_j(x,y)=\delta_{x,y}, \quad y\in V,
\]
which implies that
$G_j(x,\cdot)\in H^s(\tilde U,\Lambda^k T^* \tilde U)$. 
To estimate $\|G_j(x,\cdot)\|_{H^s(\tilde U,\Lambda^k T^* \tilde U)}$, let us take a properly supported parametrix $E_y$ of $\Delta_y^{(k)}$ in $V$ so that 
\[
E_y\Delta_y^{(k)}-1=R_y, 
\]
where $R_y$ is smoothing in $V$. 
Hence, 
\[
\|G_j(x,\cdot)\|_{H^s(\tilde U,\Lambda^k T^* \tilde U)}\le \|E_y\delta_{x,y}\|_{H^s(\tilde U,\Lambda^k T^*\tilde U)}+\|R_yG_j(x,\cdot)\|_{H^s(\tilde U,\Lambda^kT^* \tilde U)}.
\]
As
\[
E_y:H^{s-2}_{comp}(V,\Lambda^k T^* V)\to H^s_{comp}(V,\Lambda^k T^*V)
\]
is bounded, 
we get
\[
\|E_y\delta_{x,y}\|_{H^s(\tilde U,\Lambda^k T^*\tilde U)}\le C\|\delta_{x,y}\|_{H^{s-2}(V,\Lambda^k T^*V)}.
\]
The asymptotic behavior of the Green's forms \eqref{eq_asymptotic} implies that 
\[
G_j(x,\cdot)\in C^1(V,L^1_{loc}(V,\Lambda^k T^*V)).
\]
This and the fact that $R_y$ is smoothing yield that 
\[
R_yG(x,\cdot)\in C^1(V,L^2(\tilde U,\Lambda^k T^*\tilde U)). 
\]
Hence,
\[
\|R_yG_j(x,\cdot)\|_{H^s(\tilde U,\Lambda^k T^*\tilde U)}\le \|R_yG_j(x,\cdot)\|_{L^2(\tilde U,\Lambda^k T^*\tilde U)}.
\]
This shows that the function
\begin{equation}
\label{eq_G_2}
V\ni x\mapsto G_j(x,\cdot)\in H^s(\tilde U,\Lambda^k T^*\tilde U)
\end{equation}
is $C^1$. 
Combining \eqref{eq_G_1} and \eqref{eq_G_2}, we get the claim.  

\end{proof}

\begin{lem} 
\label{lem_embedding_main}
The map $\mathcal{G}_j:\tilde M_j\to H^s(\tilde U,\Lambda^k T^* \tilde U)$ defined by \eqref{eq_embedding} is an embedding. Moreover, $\mathcal{G}_j$ is  real-analytic on $\tilde M_j\setminus\overline{\tilde U}$.  

\end{lem}

\begin{proof} As the manifold $\tilde M_j$ is compact,  the fact that $\mathcal{G}_j$ is an embedding is equivalent to the fact that $\mathcal{G}_j$ is an injective immersion. 

Let us first show that $\mathcal{G}_j$ is immersion.  Suppose that there is a point $x_0\in \tilde M_j$ such that the derivative of $\mathcal{G}_j$,
\begin{align*}
&D\mathcal{G}_j(x_0):T_{x_0}\tilde M_j\to H^s(\tilde U,\Lambda^k T^* \tilde U),\\
&D\mathcal{G}_j(x_0)v=v^i\frac{\p}{\p x^i}G_j(x,\cdot)|_{x_0},
\end{align*}
where $v=v^i(\p/\p x^i)\in T_{x_0}\tilde M_j$, is not injective. Thus, there is a vector  $0\ne v\in T_{x_0}\tilde M_j$ such that
$v^i(\p/\p x^i)G_j(x_0,y)=0$ for all $y\in\tilde  U$. By the real-analyticity of the Green's forms, we have that 
$v^i(\p/\p x^i)G_j(x_0,y)=0$ for all $y\in \tilde M_j\setminus\{x_0\}$.
Since this contradicts  the asymptotic relation \eqref{eq_asymptotic} when $y$ tends to $x_0$, we get that $\mathcal{G}_j$ is immersion.

Let us now show that  $\mathcal{G}_j$ is injective. Assume that  this is not the case, then there are $x_1\ne x_2$ in $\tilde M_j$ such that 
\begin{equation}
\label{eq_injectivity}
G_j(x_1,y)= G_j(x_2,y)
\end{equation}
 for all $y\in \tilde U$. By analyticity, \eqref{eq_injectivity} holds for all $y\in \tilde M_j\setminus\{x_1,x_2\}$.
Now the asymptotic \eqref{eq_asymptotic} implies that $G_j(x_1,\cdot)$ is singular only at $y=x_1$ and $G_j(x_2,\cdot)$ is singular only at $y=x_2$. Thus, $x_1=x_2$.

\end{proof}

In the rest of this section we shall prove the following result which implies Theorem  \ref{thm_main}.

\begin{thm}
\label{thm_main_2}
 Assume that the Green's forms $G_j(x,y)$ satisfy
\begin{equation}
\label{eq_green_on_U}
G_1(x,y)=G_2(x,y), \quad (x,y)\in \tilde U\times \tilde U\setminus\{x=y\},
\end{equation}
where $\tilde U\subset \subset U$.
Thus, 
\[
\mathcal{G}_1(\tilde M_1)=\mathcal{G}_2(\tilde M_2)\subset H^s(\tilde U,\Lambda^k T^* \tilde U)
\]
and the map 
\[
(\mathcal{G}_2)^{-1}\mathcal{G}_1: \tilde M_1\to \tilde M_2
\]
is an isometry. 
\end{thm}

Before proving  Theorem \ref{thm_main_2} let us introduce some notation. Let 
\begin{align*}
N(\varepsilon_0)&=\{x\in \tilde M_1:d_{\tilde M_1}(x,\p \tilde M_1)\le \varepsilon_0\},\\
C(\varepsilon_0)&=\{x\in \tilde M_1:d_{\tilde M_1}(x,\p \tilde M_1)>\varepsilon_0\},
\end{align*}
where $\varepsilon_0>0$ is arbitrary and small enough so that $C(\varepsilon_0)$ is connected. Let $x_0\in \tilde U\cap C(\varepsilon_0)$ and $B_1\subset C(\varepsilon_0)$ be the largest connected open set containing $x_0$ such that $\mathcal{G}_1(x)\in \mathcal{G}_2(\tilde M_2)$ for all $x\in B_1$. 
The existence of the set $B_1$ follows from
\eqref{eq_green_on_U}.
Thus, as $\mathcal{G}_2$ is injective, we can define the map
\[
J=(\mathcal{G}_2)^{-1}\mathcal{G}_1: B_1\to \tilde M_2.
\]
Now let $D_1\subset B_1$ be the largest connected open set containing $x_0$ for which $J$ is a local isometry, i.e., $g_1=J^*g_2$ and $J$ is a real-analytic local diffeomorphism. 
The existence of the set $D_1$ follows from the construction of the manifold $U$.  

In order to prove Theorem \ref{thm_main_2} we will show that $ N(\varepsilon_0)\cup D_1=\tilde M_1$.  Assume the contrary, i.e. $\tilde M_1\setminus (N(\varepsilon_0)\cup D_1)\ne\emptyset$. Let 
$x_1\in \tilde M_1\setminus (N(\varepsilon_0)\cup D_1)$ be a point closest to $x_0$. Clearly $x_1\in \p D_1$.

\begin{lem}
\label{lem_limit}
There exists $x_2\in \tilde M_2^{int}$ such that 
\[
\mathcal{G}_2(x_2)=\mathcal{G}_1(x_1)=u\in H^s(\tilde U,\Lambda^k T^* \tilde U),
\]
and there exists a sequence $p_i\in D_1$ such that 
\[
\lim_{i\to\infty}p_i=x_1, \quad \lim_{i\to \infty} J(p_i)=x_2.
\]
 
\end{lem}

\begin{proof}
As $x_1\in\p D_1$, there is a sequence $p_i\in D_1$ such that $p_i\to x_1$ when $i\to\infty$. Since $J$ is a local isometry on $D_1$, there is a sequence $q_i\in \tilde M_2$ such that $\mathcal{G}_2(q_i)=\mathcal{G}_1(p_i)$. 
Since $\tilde M_2$ is compact, the sequence $q_i$ has a convergent subsequence. If there is a convergent subsequence of $q_i$ which converges to an interior point of $\tilde M_2$, then by continuity of $\mathcal{G}_1$ and $\mathcal{G}_2$, we get the claim of the lemma. 

Assume now that for all convergent subsequences of $q_i$, we have
\begin{equation}
\label{eq_q_0}
q_{i_k}\to q_0\in \p \tilde M_2, \quad k\to \infty.
\end{equation}
Lemma \ref{lem_smooth} implies that for a fixed point $y\in \tilde U$, 
\[
G_2(q_{i_k},y)\to G_2(q_0,y)=0, \quad k\to\infty.
\]
As $G_1(p_{i_k},y)=G_2(q_{i_k},y)$ for $y\in \tilde U$, we have 
 $G_1(x_1,y)=0$ for all $y\in \tilde U$. By analyticity, $G_1(x_1,y)=0$ for  all $y\ne x_1$. 
But this contradicts  the asymptotic behavior  \eqref{eq_asymptotic} of the Green's form near $y=x_1$. Hence,  
\eqref{eq_q_0} cannot be valid.

\end{proof}

Now we assume that $\tilde U\subset\subset U$ is chosen so that  
$x_j\not\in\overline{\tilde U}$, $j=1,2$. Then $\mathcal{G}_j$ is an analytic embedding in a neighborhood of $x_j$. Notice that $\tilde U\subset D_1$.

\begin{lem} We have
\[
D\mathcal{G}_1(x_1)(T_{x_1}\tilde M_1)=D\mathcal{G}_2(x_2)(T_{x_2}\tilde M_2)\subset H^s(\tilde U,\Lambda^k T^*\tilde U).
\] 
\end{lem}

\begin{proof}
By Lemma \ref{lem_limit}, there is a sequence $(p_i)_{i=1}^\infty$, $p_i\in D_1$, such that $p_i\to x_1$ when $i\to \infty$ and $J(p_i)\to x_2$ when $i\to\infty$. 
 Now by the definition of the set $D_1$, the maps 
$\mathcal{G}_1$ and $\mathcal{G}_2J$ coincide in $D_1$ and hence so do  their differentials, i.e.
\[
D\mathcal{G}_1|_{D_1}=D\mathcal{G}_2|_{J(D_1)}.
\]
Thus,
\[
D\mathcal{G}_1(p_j)(T_{p_i}\tilde M_1)=D\mathcal{G}_2(J(p_j))(T_{J(p_i)}\tilde M_2).
\]
Since by the choice of $s$, the differentials $D\mathcal{G}_1$ and $D\mathcal{G}_2$ are continuous on $\tilde M_1$ and $\tilde M_2$, we get
\[
D\mathcal{G}_1(x_1)(T_{x_1}\tilde M_1)=D\mathcal{G}_2(x_2)(T_{x_2}\tilde M_2).
\]

\end{proof}

We set
\[
\mathcal{V}=D\mathcal{G}_1(x_1)(T_{x_1}\tilde M_1)=D\mathcal{G}_2(x_2)(T_{x_2}\tilde M_2)\subset H^s(\tilde U,\Lambda^k T^*\tilde U).
\]
Since $D\mathcal{G}_j$ are injective, we have
\[
\dim \mathcal{V}= \dim T_{x_j}\tilde M_j=n. 
\]

Let 
\[
P:H^s(\tilde U,\Lambda^k T^*\tilde U)\to \mathcal{V}
\]
be the orthogonal  projection, with respect to the Hilbert space structure of $H^s$. Consider the map
\[
P\mathcal{G}_j:\tilde M_j\to \mathcal{V}, \quad x\mapsto PG_j(x,\cdot).
\]
The derivative of this map at the point $x_j$ is 
\[
DP\mathcal{G}_j(x_j)=P(D\mathcal{G}_j(x_j))=D\mathcal{G}_j(x_j),
\]
which is bijective. 
By the inverse function theorem, 
\[
P\mathcal{G}_j:\textrm{neigh}(x_j,\tilde M_j)\to \textrm{neigh}(u,\mathcal{V})
\]
is a real-analytic diffeomorphism. Writing 
\[
\mathcal{G}_j(x)=P\mathcal{G}_j(x)+(1-P)\mathcal{G}_j(x),
\]
we see that the set $\mathcal{G}_j(\tilde M_j)$ can be represented as follows,
\[
(P\mathcal{G}_j(\tilde M_j), (1-P)\mathcal{G}_j(\tilde M_j)).
\]
Thus, locally the set
\[
\{\mathcal{G}_j(x):x\in \textrm{neigh}(x_j,\tilde M_j)\}
\]
can be represented as a graph
\[
\{(v, (1-P)\mathcal{G}_j(P\mathcal{G}_j)^{-1}v):v\in \textrm{neigh}(u,\mathcal{V})\}
\]
of the real-analytic function $ (1-P)\mathcal{G}_j(P\mathcal{G}_j)^{-1}$ in $\textrm{neigh}(u,\mathcal{V})$.

Let us now show that  $x_1$ is an interior point of $B_1$. To this end, we shall prove that for any point $\tilde x_1\in  \textrm{neigh}(x_1,\tilde M_1)$, there is a unique point $\tilde x_2\in  \textrm{neigh}(x_2,\tilde M_2)$
such that 
\begin{equation}
\label{eq_B1}
\mathcal{G}_1(\tilde x_1)=\mathcal{G}_2(\tilde x_2).
\end{equation}
As the sets $\mathcal{G}_j(\textrm{neigh}(x_j,\tilde M_j))$ are the graphs of  analytic functions, \eqref{eq_B1}  is equivalent to the fact that
$P\mathcal{G}_2(\tilde x_2)=P\mathcal{G}_1(\tilde x_1)$.
Thus, 
\[
\tilde x_2=(P\mathcal{G}_2)^{-1}(P\mathcal{G}_1)(\tilde x_1)
\]
 and \eqref{eq_B1}  follows.

Therefore, in $\textrm{neigh}(x_1,\tilde M_1)$, the map 
\begin{equation}
\label{eq_J}
J=(P\mathcal{G}_2)^{-1}(P\mathcal{G}_1)
\end{equation}
is real-analytic.

To summarize, we have proved that the map
\[
J:D_1\cup \textrm{neigh}(x_1,\tilde M_1)\to \tilde M_1
\]
is real-analytic. Here we assume that the set $\textrm{neigh}(x_1,\tilde M_1)$ is connected. Moreover, when restricted  to $\tilde U\subset D_1$, 
$J$ is a local isometry. An application of 
\cite[Lemma 3', p. 256]{KobNom_book} allows us to conclude that $J$ is a local isometry on $D_1\cup \textrm{neigh}(x_1,\tilde M_1)$. Therefore, $x_1$ is an interior point of $D_1$.  This contradicts our assumption and, thus,
$N(\varepsilon_0)\cup D_1=\tilde M_1$. As $\varepsilon_0$ can be chosen arbitrarily small and by the choice of $D_1$, we have that $\mathcal{G}_1(\tilde M_1)=\mathcal{G}_2(\tilde M_2)$. Thus, $J:\tilde M_1\to \tilde M_2$ is bijective, and, therefore,  $J$ is an isometry. This proves  
Theorem \ref{thm_main_2} and, hence, 
Theorem \ref{thm_main}.

\textbf{Remark.} 
The point of the following remark is to provide an alternative approach to showing that $J$ is an isometry in a neighborhood of $x_1$, by relying upon the asymptotic behavior of the Green's forms. 
Here we follow an idea of \cite{LTU03}. 

Let us observe that \eqref{eq_B1} and \eqref{eq_J} give that 
\[
G_2(J(x),J(y))=G_1(x,y), \quad x\in \textrm{neigh}(x_1,\tilde M_1),\ y\in \tilde U.
\]
Here we used that $J(y)=y$ in $\tilde U$. As $J$ is real-analytic in $D_1\cup \textrm{neigh}(x_1,\tilde M_1)$ and Green's forms are real-analytic, we get
\begin{equation}
\label{eq_GJ}
G_2(J(x),J(y))=G_1(x,y), \quad x,y\in \textrm{neigh}(x_1,\tilde M_1),\ x\ne y.
\end{equation}

Using the asymptotic behavior  \eqref{eq_asymptotic}  of the Green's forms when $x$ is near $y$,  we shall prove that  $J$ is an isometry in the neighborhood of $y$. 
In order to do that let us recall the relation between the Riemannian metric tensor and distance function. 
For any $\xi\in T_{y}\tilde M_1$, there is a unique geodesic $\mu :\R\to \tilde M_1$ such that $\mu(0)=y$, $\dot{\mu}(0)=\xi$, and 
\[
d_{\tilde M_1}(\mu(0),\mu(t))=\int_{0}^t\sqrt{g_1(\dot{\mu}(s),\dot{\mu}(s))}ds,\quad \textrm{for } |t| \textrm{ small enough}.
\]
Thus,
\begin{equation}
\label{eq_metr}
\lim_{t\to 0}\frac{d_{\tilde M_1}(\mu(0),\mu(t))}{t}=g_1(\xi,\xi)^{1/2}.
\end{equation}

It follows from  \eqref{eq_asymptotic} that as $x\to y$, we have 
\begin{align*}
 G_2(J(x),J(y)) \sim&\frac{1}{(n-2)s_nd_{\tilde M_2}(J(x),J(y))^{n-2}} \Gamma_{(i_1,\dots,i_k),(j_1,\dots,j_k)}^{g_2}(J(y))\\
 & (L_{(i_1,\dots,i_k),(j_1,\dots,j_k)}(y))^2dx^{i_1}\wedge \dots \wedge dx^{i_k}\cdot dy^{j_1}\wedge \dots \wedge dy^{j_k},
\end{align*}
where $L_{(i_1,\dots,i_k),(j_1,\dots,j_k)}(y)$ are certain $k\times k$-minors of the Jacobian of $J(y)$.  There are indices  $(i_1,\dots,i_k),(j_1,\dots,j_k)$ such that 
$L_{(i_1,\dots,i_k),(j_1,\dots,j_k)}(x_1)\ne 0$. Thus, shrinking the neighborhood of $x_1$, if necessary, we may assume that
\[
L_{(i_1,\dots,i_k),(j_1,\dots,j_k)}(y)\ne 0, \quad \forall y\in \textrm{neigh}(x_1,\tilde M_1).
\]
Hence, \eqref{eq_GJ} implies that
\begin{equation}
\label{eq_metr_2}
\lim_{x\to y}\frac{d_{\tilde M_2}(J(x),J(y))}{d_{\tilde M_1}(x,y)}=f(y),
\end{equation}
where 
\[
f(y)=
\left(\frac{\Gamma_{(i_1,\dots,i_k),(j_1,\dots,j_k)}^{g_2}(J(y))
 (L_{(i_1,\dots,i_k),(j_1,\dots,j_k)}(y))^2}{\Gamma_{(i_1,\dots,i_k),(j_1,\dots,j_k)}^{g_1}(y)}\right)^{1/(n-2)}
\]
 is positive real-analytic function in $\textrm{neigh}(x_1,\tilde M_1)$.  

Combining \eqref{eq_metr} and \eqref{eq_metr_2}, we get
\begin{equation}
\label{eq_metr_3}
\frac{g_2(J(y))(J'(y)\xi,J'(y)\xi)}{g_1(y)(\xi,\xi)}=f(y)^2, 
\end{equation}
for any $y\in \textrm{neigh}(x_1,\tilde M_1)$ and $\xi\in T_{y}\tilde M_1$. 
As both right and left-hand sides of \eqref{eq_metr_3} are real-analytic in $\textrm{neigh}(x_1,\tilde M_1)\cup D_1$, \eqref{eq_metr_3}  holds also on $D_1$. 
Since
\[
J^*g_2=g_1, \quad\textrm{on}\quad \tilde U\subset D_1,
\]
we have $f(y)\equiv 1$ on $\textrm{neigh}(x_1,\tilde M_1)$. 
Hence, $J$ is an isometry in a neighborhood of $x_1$.

\section{The case of a compact $2$-dimensional manifold. Proof of Theorem \ref{thm_main_2dim}}

\label{sec_2dim}

Let $(M,g)$ be a smooth compact connected Riemannian manifold of dimension $n=2$ with a smooth non-empty boundary.  
 Let $(x^1,x^2)$ be the boundary normal coordinates defined locally near a point at the boundary. Here $x^1$ is a local coordinate for $\p M$ and $x^2\ge 0$ is the distance to the boundary. In these coordinates, the metric has the following form
\[
g=g_{11}(x^1,x^2)(dx^1)^2+(dx^2)^2.
\]
Recall that the Hodge star isomorphism is defined by
 \begin{equation}
 \label{eq_hodge_star}
 \omega\wedge *\omega=g(\omega,\omega)\mu,
 \end{equation}
where $\mu\in C^\infty (M,\Lambda^2 T^*M)$ is the volume form given in the boundary normal coordinates by
\[
\mu=\sqrt{g_{11}}dx^1\wedge dx^2,
\]
assuming that $dx^1$, $dx^2$ is a positive basis of $T^*_x M$.
Thus, \eqref{eq_hodge_star} implies that
\[
dx^1\wedge * dx^1=g^{11}\sqrt{g_{11}}dx^1\wedge dx^2,
\]
and therefore,
\begin{equation}
\label{eq_hodge_star2}
*dx^1=\sqrt{g^{11}}dx^2,\quad *dx^2=-\sqrt{g_{11}}dx^1, \quad *(dx^1\wedge dx^2)=\sqrt{g^{11}}. 
\end{equation}
Here $g^{11}=g_{11}^{-1}$.

\subsection{The case of $0$-forms} Let  $\tilde g$ be a smooth Riemannian metric on $M$ in the same conformal class as $g$, i.e. there is a smooth real-valued function $\varphi$ on $M$ such that 
\[
  \tilde g = e^{2\varphi}g.
\]
Then it is known \cite[Chapter 1.J]{Bes_book} that under the conformal change, the Hodge star operator on $k$-forms satisfies 
\begin{equation}
\label{eq_comf_1}
    *_{\tilde g} = e^{(2-2k)\varphi}*_g.
\end{equation}

It follows from \eqref{eq_comf_1} that the Hodge Laplacian on $0$-forms is conformally invariant, i.e.
\[
\Delta^{(0)}_{\tilde g}=e^{-2\varphi}\Delta^{(0)}_{g}.
\]
Moreover, $\mathcal{C}^{(0)}_{\tilde g}=\mathcal{C}^{(0)}_{g}$. 

When computing the set of the Cauchy data $\mathcal{C}^{(0)}_{g}$ in boundary normal coordinates, we get, when $\omega\in C^\infty(M)$,
\begin{align*}
\mathbf{t}\omega&=\omega|_{x^2=0},\\
\mathbf{n}_gd\omega&=-\sqrt{g_{11}(x^1,0)}\p_{x^2}\omega|_{x^2=0} dx^1. 
\end{align*}
Notice that the set of the Cauchy data $\mathcal{C}^{(0)}_{g}$ is the graph of the Dirichlet-to-Neumann map as defined in \cite{LeeUhl89}.   
As explained in  \cite{LeeUhl89},  $\mathcal{C}^{(0)}_{g}$ does not contain any information about the metric along the boundary. 
However, if one considers the Dirichlet-to-Neumann map as a function valued map by setting
\[
\omega|_{x^2=0}\mapsto \p_{x^2}\omega|_{x^2=0},
\]
then it contains the information about the restriction of the metric to the boundary, see \cite{LeeUhl89}.

\subsection{The case of $1$-forms}
Given a $1$-form $\omega=\omega_1dx^1+\omega_2dx^2$, we shall identify $\omega$ with the vector $(\omega_1,\omega_2)^{\mathrm{t}}$. An explicit computation using \eqref{eq_hodge_star2} shows that in boundary normal coordinates, the   Hodge Laplacian on $1$-forms has the form
\begin{equation}
\label{eq_lap_1}
\Delta^{(1)}=D_{x^2}^2I + g^{11}D_{x^1}^2I +iE(x)D_{x^2}+i F(x)D_{x^1}+ Q(x),
\end{equation}
where $D_{x^j}=\frac{1}{i}\p_{x^j}$, $j=1,2$, and
\[
E(x)=\frac{\p _{x^2} g^{11}}{2g^{11}}\begin{pmatrix} -1 & 0\\ 0& 1\end{pmatrix},
\]
\[
F(x)=\begin{pmatrix}
-\frac{3}{2}\p _{x^1}g^{11} & \frac{\p_{x^2}g^{11}}{g^{11}}\\
-\p_{x^2}g^{11}& -\frac{\p _{x^1}g^{11}}{2}
\end{pmatrix},
\quad
Q(x)=\begin{pmatrix} -\frac{\p_{x^1}^2 g^{11}}{2} & \frac{1}{2}\p_{x^1x^2}^2\log g^{11},\\
-\frac{\p^2_{x^1x^2} g^{11}}{2} & \frac{1}{2}\p^2_{x^2}\log g^{11}
\end{pmatrix}.
\]
Here $I$ is the $2\times 2$ identity matrix.

Our next goal is to show that the knowledge of the Cauchy data for harmonic $1$-forms determines the Taylor series at the boundary of the metric $g$ in the boundary normal coordinates.
Let us point out that the reconstruction of the Taylor series in this case becomes possible, thanks to the special structure of the lower order terms in $\Delta^{(1)}$. 
We shall need the following lemma.

\begin{lem}
\label{lem_factor}
There exists a matrix-valued  pseudodifferential operator $A(x,D_{x^1})$ of order one in $x^1$ depending smoothly on $x^2$ such that 
\begin{equation}
\label{eq_lap_2}
\Delta^{(1)}=(D_{x^2}I+iE(x)-iA(x,D_{x^1}))(D_{x^2}I+iA(x,D_{x^1})),
\end{equation}
modulo a smoothing operator. Here $A(x,D_{x^1})$ is unique modulo a  smoothing term, if we require that its principal symbol satisfies 
\[
\sigma(A((x^1,0),D_{x^1}))=-\sqrt{g^{11}(x^1,0)}|\xi_1| I.
\]
\end{lem} 

\begin{proof}
The existence of the factorization \eqref{eq_lap_2}  follows closely  \cite[Proposition 1.1]{LeeUhl89}, where the case of the Laplacian on functions is considered. 

Combining \eqref{eq_lap_1} and \eqref{eq_lap_2}, we see that
\begin{equation}
\label{eq_lap_3}
A^2+i[D_{x^2}I ,A]-E(x)A=g^{11}D_{x^1}^2I+iF(x)D_{x^1}+Q(x),
\end{equation}
modulo a smoothing operator. Let us write the full symbol of $A(x,D_{x^1})$ as follows,
\[
A(x,\xi_1)\sim \sum_{j=-\infty}^1 a_j(x,\xi_1),
\]
with $a_j$ taking values in $2\times 2$ matrices with entries homogeneous of degree $j$ in $\xi_1$. Thus, 
\eqref{eq_lap_3} implies that
\begin{align}
\label{eq_lap_4}
\sum_{l=-\infty}^2\underset{\alpha\ge 0,\ j,k\le 1}{(\sum_{j+k-\alpha=l}}\frac{1}{\alpha!}\p_{\xi_1}^\alpha a_jD_{x^1}^\alpha a_k)
-E(x)\sum_{j=-\infty}^1 a_j(x,\xi_1)+\sum_{j=-\infty}^1 \p_{x^2}a_j(x,\xi_1)\\ \nonumber
=g^{11}\xi_1^2I+iF(x)\xi_1+Q(x).
\end{align}
Equating the terms homogeneous of degree two in \eqref{eq_lap_4}, we get
\[
a_1^2(x,\xi_1)=g^{11}(x)\xi_1^2I,
\] 
so we should choose $a_1(x,\xi_1)$ to be the scalar matrix given by  
\[
a_1(x,\xi_1)=-\sqrt{g^{11}(x)}|\xi_1|I.
\]
Equating the terms homogeneous of degree one in \eqref{eq_lap_4}, we have an equation for $a_0(x,\xi_1)$,
\begin{equation}
\label{eq_a_0}
2a_1a_0=iF(x)\xi_1+E(x)a_1-\p_{x^2}a_1 -\p_{\xi_1}a_1D_{x^1}a_1,
\end{equation}
which is uniquely solvable for $a_0$ given our choice of $a_1$.  We shall now consider terms homogeneous of degree zero in \eqref{eq_lap_4},   
\begin{equation}
\label{eq_a_-1}
2a_1a_{-1}=Q(x)+E(x)a_0-\p_{x^2}a_0-\underset{\alpha\ge 0,\ 0\le  j,k\le 1}{\sum_{j+k-\alpha=0}}\frac{1}{\alpha!}\p_{\xi_1}^\alpha a_jD_{x^1}^\alpha a_k),
\end{equation}
which is uniquely solvable for $a_{-1}$  given our choice of $a_1$ and $a_0$.

Proceeding recursively with respect to the degree of homogeneity in \eqref{eq_lap_4}, we choose, for $m\ge 1$, 
\[
a_{-m-1}=\frac{1}{2}a_1^{-1}(Ea_{-m}-\p_{x^2}a_{-m}-\underset{\alpha\ge 0,-m\le j,k\le 1}{\sum_{j+k-\alpha=-m}}\frac{1}{\alpha!}\p_{\xi_1}^\alpha a_jD_{x^1}^\alpha a_k).
\]
This completes the proof.
\end{proof}

\begin{prop}
\label{prop_metric}
The knowledge of the subset of the set of Cauchy data on $1$-forms given by
\[
\{(\mathbf{t}\omega,\mathbf{n}d\omega):\omega\in C^\infty(M, \Lambda^1 T^* M), \Delta^{(1)}\omega=0, \mathbf{n}\omega=0\},
\]
determines the full Taylor series at the boundary of the metric $g$ in the boundary normal coordinates.

\end{prop}

\begin{proof}

We shall first show how to recover the restriction of the metric $g$ to the boundary. To this end, we shall compute the Cauchy data of  a harmonic $1$-form $\omega=\omega_1dx^1+\omega_2dx^2$
 in boundary normal coordinates. Using \eqref{eq_hodge_star2}, we get
\begin{align*}
\mathbf{t}\omega&=\omega_{1}(x_1,0)dx^1,\\
\mathbf{n}\omega&=-\omega_2(x_1,0)\sqrt{g_{11}(x^1,0)}dx^1,\\
\mathbf{n}d\omega&=\left(\frac{\p \omega_2}{\p x^1}-\frac{\p \omega_1}{\p x^2}\right)|_{x^2=0}\sqrt{g^{11}(x^1,0)},\\
\mathbf{t}\delta\omega&=-\left(\frac{\p}{\p x^1}(\omega_1\sqrt{g^{11}})+ \frac{\p}{\p x^2}(\omega_2\sqrt{g_{11}})\right)|_{x^2=0}\sqrt{g^{11}(x^1,0)}.
\end{align*}

The same argument as in \cite[Proposition 1.2]{LeeUhl89} and in \cite{JL05} shows that the operator $A(x,D_{x^1})$ in  the factorization \eqref{eq_lap_2} of the Hodge Laplacian $\Delta^{(1)}$ has the following meaning: given a harmonic $1$-form $\omega=\omega_1dx^1+ \omega_2dx^2$, we have
\[
\begin{pmatrix}
\p_{x^2}\omega_1|_{x^2=0}\\
\p_{x^2}\omega_2|_{x^2=0}
\end{pmatrix}=A((x^1,0),D_{x^1})
\begin{pmatrix}
\omega_1|_{x^2=0}\\
\omega_2|_{x^2=0}
\end{pmatrix},
\]
modulo a smoothing  operator.   
We shall represent the operator $A$ by a $2\times 2$-matrix of scalar operators $A_{jk}$, $1\le j,k\le 2$. 
Taking an arbitrary harmonic $1$-form $\omega$ with  $\mathbf{n}\omega=0$, we get
\[
\p_{x^2}\omega_1|_{x^2=0}=A_{11}((x^1,0),D_{x^1})\omega_1|_{x^2=0},
\]
modulo smoothing. Thus,
\[
\sqrt{g^{11}(x^1,0)}\p_{x^2}\omega_1|_{x_2=0}=\sqrt{g^{11}(x^1,0)}A_{11}((x^1,0),D_{x^1})\omega_1|_{x^2=0},
\]
modulo smoothing, and, therefore, the knowledge of the Cauchy data for harmonic $1$-forms $\omega$ with  $\mathbf{n}\omega=0$ implies the knowledge of the full symbol of the operator
\begin{equation}
\label{eq_metric_1}
\sqrt{g^{11}(x^1,0)}A_{11}((x^1,0),D_{x^1}).
\end{equation}
The factorization \eqref{eq_lap_2} yields that the principal symbol of $A_{11}((x^1,0),D_{x^1})$ is $-\sqrt{g^{11}(x^1,0)}|\xi_1|$. Hence, from the knowledge of the principal symbol of the operator \eqref{eq_metric_1}, we determine the metric $g_{11}(x^1,0)$ along the boundary. Having found the metric along the boundary, we recover the full symbol of $A_{11}((x^1,0),D_{x^1})$.

We shall determine the normal derivatives of the metric at the boundary from the knowledge of the full symbol of the operator $A_{11}((x^1,0),D_{x^1})$. 
It follows from \eqref{eq_a_0} that the homogeneous symbol of order zero $a_{0,11}$ of the operator $A_{11}$ is given by
\begin{equation}
\label{eq_a_0_2}
a_{0,11}=-\frac{1}{2\sqrt{g^{11}}|\xi_1|}((Ea_1)_{11}+\p_{x^2}\sqrt{g^{11}}|\xi_1|)+T_0=
-\frac{\p_{x^2}g^{11}}{2g^{11}}+T_0,
\end{equation}
where 
\[
T_0=-\frac{1}{2\sqrt{g^{11}}|\xi_1|}(iF_{11}\xi_1-(\p_{\xi_1}a_1D_{x^1}a_1)_{11})
\]
contains only the metric $g_{11}(x^1,0)$ and its tangential derivatives. 
Thus, using \eqref{eq_a_0_2}, we can determine the normal derivative $\p_{x^2} g_{11}(x^1,0)$.  Notice that \eqref{eq_a_0} implies that $a_0((x^1,0),\xi_1)$ contains only the metric along the boundary $g_{11}(x^1,0)$, its tangential derivatives  and the first normal derivative $\p_{x^2}g_{11}(x^1,0)$. 

It follows from \eqref{eq_a_-1} that the homogeneous symbol of order $-1$ of the operator $A_{11}$ has the form
\begin{equation}
\label{eq_a_-1_2}
a_{-1,11}=-\frac{1}{2\sqrt{g^{11}}|\xi_1|}(-\p_{x^2}a_{0,11})+T_{-1},
\end{equation}
where 
\[
T_{-1}=-\frac{1}{2\sqrt{g^{11}}|\xi_1|}(Q_{11}+(Ea_0)_{11}-\underset{\alpha\ge 0,\ 0\le  j,k\le 1}{\sum_{j+k-\alpha=-0}}\frac{1}{\alpha!}(\p_{\xi_1}^\alpha a_jD_{x^1}^\alpha a_k)_{11})
\]
contains only the metric along the boundary, its tangential derivatives  and its first normal derivative $\p_{x^2}g_{11}(x^1,0)$. 
Substituting $a_{0,11}$ into \eqref{eq_a_-1_2}, we get
\begin{equation}
\label{eq_a_-1_3}
a_{-1,11}=-\frac{1}{2\sqrt{g^{11}}|\xi_1|}\left(\frac{\p^2_{x^2}g^{11}}{2g^{11}}\right)+\tilde T_{-1},
\end{equation}
where $\tilde T_{-1}$ contains only the metric along the boundary, its tangential derivatives  and its first normal derivative. 
Hence, from \eqref{eq_a_-1_3}, we can recover the second normal derivative of the metric along the boundary.
Notice that $a_{-1}$ contains only the metric, its tangential derivatives and its normal derivatives up to order $2$ along the boundary.

Let us now proceed by induction with respect to the degree of homogeneity in \eqref{eq_lap_4}. Let $m\ge 1$ and suppose we have shown that
\begin{equation}
\label{eq_lap_6}
a_{-j,11}=-\frac{1}{(2\sqrt{g^{11}}|\xi_1|)^j}\left( \frac{\p_{x^2}^{j+1}g^{11}}{2g^{11}}\right)
+T_{-j}, \quad 1\le  j\le m,
\end{equation}
where $T_{-j}$ is an expression involving only  $\p_{x^2}^{l} g_{11}$, $l=0,\dots, j$, and their tangential derivatives along $\p M$, and,
furthermore, also suppose that we have shown that $a_{-j}$ contains only $\p_{x^2}^{l} g_{11}$, $l=0,\dots, j+1$, and their tangential derivatives along $\p M$.

We shall show how to determine the normal derivative $\p_{x^2}^{m+2} g_{11}(x^1,0)$ 
from the knowledge of the term $a_{-m-1,11}$ of the full symbol of $A_{11}$. 
 We have
\begin{equation}
\label{eq_lap_7}
a_{-m-1,11}=-\frac{1}{2\sqrt{g^{11}}|\xi_1|}((Ea_{-m})_{11}-\p_{x^2}a_{-m,11}-\underset{\alpha\ge 0,-m\le j,k\le 1}{\sum_{j+k-\alpha=-m}}\frac{1}{\alpha!}(\p_{\xi_1}^\alpha a_jD_{x^1}^\alpha a_k)_{11})
\end{equation}
It follows from \eqref{eq_lap_6} that the only term in the right hand side of \eqref{eq_lap_7} which contains the normal derivative $\p_{x^2}^{m+2} g_{11}(x^1,0)$ is 
\[
\frac{1}{2\sqrt{g^{11}}|\xi_1|}\p_{x^2}a_{-m,11}.
\] 
We conclude that the representation \eqref{eq_lap_6} holds also for $a_{-m-1}$ and, thus, the derivative $\p_{x^2}^{m+2} g_{11}(x^1,0)$ is determined.

\end{proof}

We are now ready to complete the \textbf{proof of  Theorem \ref{thm_main_2dim}}. To this end, let us recall the asymptotics of the Green's $1$-form $G(x,y)$, as $x\to y$, $y\in \tilde M^{int}$,
\begin{equation}
\label{eq_asyptotics_2}
G(x,y)\sim -\frac{1}{s_2} \log d_{\tilde M}(x,y) \det(g) dx^{1}\wedge dx^{2}\cdot dy^{1}\wedge dy^{2},
\end{equation}
see  \cite[p.136]{DS52}. We refer to  \eqref{eq_asymptotic} for the notation.  Given Proposition \ref{prop_metric} and the asymptotics \eqref{eq_asyptotics_2},
the proof of Theorem \ref{thm_main_2dim} is obtained by repeating the proof of Theorem \ref{thm_main}.

\subsection{The case of $2$-forms}
Let $g$ and $\tilde g$ be two smooth Riemannian metrics on $M$ in the same conformal class, i.e. there is a smooth real-valued function $\varphi$ on $M$ such that 
\[
  \tilde g = e^{2\varphi}g.
\]
Then using \eqref{eq_comf_1}, for the Hodge Laplacian on $2$-forms,  we have
\begin{align*}
&\Delta^{(2)}_{\tilde g}\omega=\Delta^{(2)}_g(e^{-2\varphi}\omega),\quad \omega\in C^\infty(M,\Lambda^2 T^*M),\\
&\mathbf{n}_{\tilde g}\omega=\mathbf{n}_{g}(e^{-2\varphi}\omega),\\
&\mathbf{t}\delta_{\tilde g}\omega=\mathbf{t}\delta_g(e^{-2\varphi}\omega).
\end{align*}
Hence, 
\begin{equation}
\label{eq_cauchy_two_new}
\mathcal{C}^{(2)}_{\tilde g}=\mathcal{C}^{(2)}_{g}.
\end{equation}

\textbf{Remark.} Let $\omega=\omega_{12}dx^1\wedge dx^2$ be a $2$-form.  An explicit computation using \eqref{eq_hodge_star2} shows that in boundary normal coordinates, the   Hodge Laplacian on $2$-forms has the form
\begin{equation}
\label{eq_lap_two_forms}
\Delta^{(2)}=D_{x^2}^2 + g^{11}D_{x^1}^2 +iE(x)D_{x^2}+i F(x)D_{x^1}+ Q(x),
\end{equation}
where 
\[
E(x)=-\frac{\p_{x^2}g^{11}}{2g^{11}},\quad F(x)=-\frac{3}{2}\p_{x^1}g^{11},\quad
Q(x)=-\frac{\p^2_{x^1}g^{11}}{2}-\frac{\p^2_{x^2}\log g^{11}}{2}.
\]
Factorizing the Laplacian $\Delta^{(2)}$ similarly to the case of $1$-forms, one can then determine the principal symbol of the corresponding Dirichlet-to-Neumann map,  and show that it does not carry any information about the metric on the boundary of the manifold.

\section{The case of a complete manifold. Proof of Theorem \ref{thm_main_complete}}

\label{sec_complete}

In order to prove Theorem \ref{thm_main_complete} we shall follow the proof of Theorem  \ref{thm_main}, valid in the compact case.
It follows from \cite{JL05}, when $n\ge 3$ and $k=1,\dots,n$, and from Proposition \ref{prop_metric}, when $n=2$ and $k=1$, that the knowledge of the Dirichlet-to-Neumann map $\Lambda^{(k)}_{g,\Gamma}$  determines all the normal derivatives $\p_{x^n}^l g_{ij}(x',0)$, $l=0,1,2,\dots$, of the metric tensor $g$ at $\Gamma$, in boundary normal coordinates.

Recall that we work under the assumption (A1), when $n\ge 3$, and the assumption (A2), when $n=2$. As in the proof of Theorem  \ref{thm_main}, we construct the real-analytic manifold $\tilde M$ by attaching to $M$ a boundary collar $\Gamma\times (-r,0]$, 
with the extension of the metric $g$. Here $r$ is small enough, so that the extended metric remains real-analytic. 
In the same way as in Proposition \ref{prop_Friedrich}, given the Hodge Laplacian on $\tilde M$, we shall  consider the non-negative self-adjoint realization $\Delta^{(k)}_{F,\tilde M}$.

Since the manifolds $\tilde M$ and $M$ differ by a compact region, the following result is essentially well known, see \cite{Bar2000}. 

\begin{prop}
\label{prop_comp_ess}

If $0\notin\textrm{spec}_{ess}(\Delta^{(k)}_{F,M})$, then 
$0\notin\textrm{spec}_{ess}(\Delta^{(k)}_{F,\tilde M})$. 
\end{prop}

Recall now from the assumptions of Theorem \ref{thm_main_complete} that  $0\notin\textrm{spec}(\Delta^{(k)}_{F,\tilde M})$. Thus, by Proposition \ref{prop_comp_ess}, we may introduce the inverse operator 
\[
(\Delta^{(k)}_{F,\tilde M })^{-1}: L^2(\tilde M,\Lambda^k T^*\tilde M)\to \mathcal{D}(\Delta^{(k)}_{F,\tilde M}),
\]
given by
\[
(\Delta^{(k)}_{F,\tilde M})^{-1} w(x)=\int_{\tilde M}G(x,y)\wedge *w(y).
\]
 Here the Schwartz kernel $G(x,y)$ is the corresponding Green's form, which satisfies 
\begin{align*}
&\Delta_{x}^{(k)}G(x,y)=\delta_{x,y}
\quad \text{in}\quad \tilde M,\\
&\mathbf{t}_xG(x,y)=0\quad \text{on}\quad \p\tilde M,\\
&\mathbf{n}_xG(x,y)=0\quad \text{on}\quad \p\tilde M,
\end{align*}
where $y\in \tilde M\setminus\p\tilde M$ and 
\[
\delta_{x,y}=\sum_{\underset{
j_1<\dots <j_k}{i_1<\dots <i_k}}\delta(y-x)dx^{i_1}\wedge \dots\wedge dx^{i_k}\cdot
dy^{j_1}\wedge \dots\wedge dy^{j_k}
\]
 is the delta double current, supported at $x=y$.

\begin{prop}

The Green's form $G(x,y)$ has the following properties. 

\begin{itemize}

\item[(i)]  $G(x,y)=G(y,x)$. 

\item[(ii)]  Let $y\in \tilde M\setminus\p \tilde M$. Then $G(x,y)$ is real-analytic for $x\in \tilde M\setminus(\{y\}\cup \p\tilde M)$. 

\item[(iii)] Let $y\in \tilde M\setminus\p \tilde M$. Then $G(x,y)$ is $C^\infty$-smooth up to the boundary of $\p \tilde M$, away from $y$. 

\item[(iv)] Let $n\ge 3$. Then as $x\to y\in \tilde M^{int}$, $G(x,y)$ has the following asymptotic behavior 
\begin{equation}
\label{eq_asymptotic_complete}
\begin{aligned}
G(x,y)\sim&\frac{1}{(n-2)s_nd_{\tilde M}(x,y)^{n-2}} \Gamma_{(i_1,\dots,i_k),(j_1,\dots,j_k)}(y)
\\
&dx^{i_1}\wedge \dots \wedge dx^{i_k}\cdot dy^{j_1}\wedge \dots \wedge dy^{j_k},
\end{aligned}
\end{equation}
where $d_{\tilde M}(x,y)$ is the geodesic distance from $x$ to $y$ in $\tilde M$, $s_n$ is the area of the unit $n$-sphere and 
 $\Gamma_{(i_1,\dots,i_k),(j_1,\dots,j_k)}(y)$ is a positive real-analytic function on $\tilde M$. The asymptotic relation \eqref{eq_asymptotic_complete} can be differentiated  with respect to $x$ any number of times.

\item[(v)] Let $n=2$. Then the asymptotics of the Green's $1$-form $G(x,y)$, as $x\to y$, $y\in \tilde M^{int}$, has the following form,
\begin{equation}
\label{eq_asyptotics_2-complete}
G(x,y)\sim -\frac{1}{s_2} \log d_{\tilde M}(x,y) \det(g) dx^{1}\wedge dx^{2}\cdot dy^{1}\wedge dy^{2}.
\end{equation}
 The asymptotic relation \eqref{eq_asyptotics_2-complete} can be differentiated  with respect to $x$ any number of times.

\end{itemize}

\end{prop}

\begin{proof}

The property (i) is clear, and (ii) follows from \cite[Theorem 4.1, p. 108]{EgoShu}, since $G(x,y)$ is a solution of  an  elliptic system of partial differential equations on $\tilde M\setminus(\{y\}\cup \p\tilde M)$ with real-analytic coefficients. 
The property (iii) follows from \cite[Theorem 3.4.6]{Sch_book}. 

Let us show (iv) and (v). Consider $\Omega\subset\subset \tilde M^{int}$. Let $G^\Omega(x,y)$ be a Green's form on a compact manifold $\overline{\Omega}$ defined as in Section \ref{sec_Green's_forms}. Then 
\[
(\Delta^{(k)}_x+\Delta^{(k)}_y)(G(x,y)-G^{\Omega}(x,y))=0, \quad (x,y)\in \Omega\times \Omega,
\]
and therefore,
\[
G(x,y)=G^\Omega(x,y)+h(x,y), \quad h(x,y)\text{ is } C^\infty-\textrm{smooth double form on}\ \Omega\times \Omega.
\]
The proof is complete. 

\end{proof}

As before, let us set $U= \tilde M_1\setminus M_1= \tilde M_2\setminus M_2$. In analogy with Lemma \ref{lem_green_new}, we have the following result.  

\begin{lem}
The Green's forms $G_j(x,y)$ satisfy
\[
G_1(x,y)=G_2(x,y), \quad (x,y)\in U\times U\setminus\{x=y\}.
\]
\end{lem}

Let  $\tilde U\subset\subset U$ be a relatively compact subset of $U$. 
In order to prove Theorem \ref{thm_main_complete}, as in the compact case, we shall make use of the maps
\begin{equation}
\label{eq_emb_complete}
\mathcal{G}_j:\tilde M_j\to H^s(\tilde U,\Lambda^k T^* \tilde U),\quad
 x\mapsto  G_j(x,\cdot),
\end{equation}
defined for some $s<1-n/2$. 
Proposition \ref{prop_regularity_emb} continues to be valid in the case when $\tilde M_j$ are complete, and therefore, 
the maps $\mathcal{G}_j$ are $ C^1(\tilde M_j,H^s(\tilde U,\Lambda^k T^* \tilde U))$, for $s<1-n/2$. Furthermore, arguing as in the proof of Lemma \ref{lem_embedding_main}, we see that the maps  \eqref{eq_emb_complete} are injective immersions, which suffices for us to proceed with the proof of Theorem \ref{thm_main_complete} along the lines of the argument presented in  the proof of Theorem \ref{thm_main}. 
The only modification required in the complete non-compact case concerns the proof of Lemma \ref{lem_limit}. We shall therefore give a proof of this result in this case.

\textbf{Proof of Lemma \ref{lem_limit} in the non-compact case}. 
As $x_1\in\p D_1$, there is a sequence $p_i\in D_1$ such that $p_i\to x_1$ when $i\to\infty$. 
Hence, 
\[
d_{\tilde M_1}(x_0,p_i)\to d_{\tilde M_1}(x_0,x_1) \quad i\to \infty.
\]
As every sequence in $\R$ has a monotone subsequence, we can assume that
\[
d_{\tilde M_1}(x_0,p_i)\le d_{\tilde M_1}(x_0,p_{i+1})\le d_{\tilde M_1}(x_0,x_1), \quad i=1,2,\dots. 
\]
Let $\mu_i$ be the shortest curve from $x_0$ to $p_i$ in $\tilde M_1$, i.e. assume that
\[
\textrm{length}_{g_1}(\mu_i)=d_{\tilde M_1}(x_0,p_i). 
\]
The fact that  $x_1$ is the closest point of $\tilde M_1\setminus(N(\varepsilon_0)\cup D_1))$ implies that $\mu_i$ is contained in $D_1$. 
Since $J$ is a local isometry on $D_1$, there is a sequence $q_i\in \tilde M_2$ such that $q_i=J(p_i)$ and, moreover,  
\[
d_{\tilde M_2}(q_i,x_0)\le \textrm{length}_{g_2}(J(\mu_i))=\textrm{length}_{g_1}(\mu_i)\le d_{\tilde M_1}(x_0,x_1). 
\]
Thus, the sequence $q_i$ is bounded and the remainder of the proof of Lemma \ref{lem_limit} proceeds similarly to the compact case. 

This completes the proof of Theorem \ref{thm_main_complete}.

\section{Acknowledgements}
The research of K.K. was financially supported by the
Academy of Finland (project 125599).  The research of M.L. 
was financially supported 
 by the Academy of Finland Center of Excellence programme 213476.
G.U. was partly supported by NSF and a Walker Family Endowed Professorship.


\begin{thebibliography} {1}

\bibitem{BelSch1008} 
Belishev, M.,  Sharafutdinov, V., \emph{Dirichlet to Neumann operator on differential forms},  Bull. Sci. Math.  \textbf{132}  (2008),  no. 2, 128--145.

\bibitem{Bes_book}
Besse, A.,  \emph{Einstein manifolds}, Reprint of the 1987 edition. Classics in Mathematics. Springer-Verlag, Berlin, 2008, 516 pp.


\bibitem{Bar2000}
B\"ar, C.,  \emph{The Dirac operator on hyperbolic manifolds of finite volume}, J. Differential Geom. \textbf{54} (2000), no. 3, 439--488. 


\bibitem{Dav_book} Davies, E. B., \emph{Spectral theory and differential operators},
   Cambridge Studies in Advanced Mathematics,
   \textbf{42}, Cambridge University Press, 1995, 182 pp.


\bibitem{DonXav84}
Donnelly, H., Xavier, F., \emph{On the differential form spectrum of negatively curved Riemannian manifolds},
Amer. J. Math. \textbf{106} (1984), no. 1, 169--185. 

\bibitem{DS52}
Duff, G. F. D., Spencer, D. C.,  \emph{Harmonic tensors on {R}iemannian manifolds with boundary},  Ann. of Math. (2)  \textbf{56} (1952), 128--156.

\bibitem{EgoShu}
Egorov, Yu. V., 
\emph{Partial differential equations. {IV},Microlocal analysis},
   Encyclopaedia of Mathematical Sciences, \textbf{33}, Springer-Verlag, Berlin, 1993, 241pp.

\bibitem{GLU03}  
Greenleaf, A., Lassas, M. and Uhlmann, G.,
Uhlmann, \emph{On
nonuniqueness for
Calder\'on's,
inverse problem},  
Math. Res. Lett., {\bf 10} (2003),
no. 5-6, 685-693.

\bibitem{GKLU09} Greenleaf, A., Kurylev, Y., Lassas, M. and Uhlmann, G.,
\emph{Invisibility and inverse problems}, Bulletin AMS, \textbf{46} (2009), no 1, 55-97. 

\bibitem{JL05}
Joshi, M. S., Lionheart, W. R. B., \emph{An inverse boundary value problem for harmonic differential forms}, Asymptot. Anal., \textbf{41} (2005), no. 2, 93--106. 

\bibitem{KobNom_book}
Kobayashi, S.,  Nomizu, K., \emph{Foundations of differential geometry}. Vol. I.  Wiley Classics Library. A Wiley-Interscience Publication. John Wiley and Sons, Inc., New York, 1996, 329 pp.


\bibitem{LTU03}
Lassas, M., Taylor, M., Uhlmann, G., \emph{The {D}irichlet-to-{N}eumann map for
complete {R}iemannian manifolds with boundary},  Comm. Anal. Geom.,  \textbf{11}  (2003),  no. 2,
207--221.

\bibitem{LU01}
Lassas, M., Uhlmann, G., \emph{On determining a {R}iemannian manifold from the
{D}irichlet-to-{N}eumann map},  Ann. Sci. \'Ecole Norm. Sup. (4),  \textbf{34}  (2001),  no. 5, 771--787.


\bibitem{LeeUhl89} 
Lee, J., Uhlmann, G., \emph{Determining anisotropic real-analytic conductivities by boundary measurements},
Comm. Pure Appl. Math. \textbf{42} (1989), no. 8, 1097--1112. 

\bibitem{Nach88}
Nachman, A.,  \emph{Reconstructions from boundary measurements},  Ann. of Math. (2) \textbf{128} (1988),
no. 3, 531--576.


\bibitem{Rham_book}
de Rham, G., \emph{Differentiable manifolds.
Forms, currents, harmonic forms}, Grundlehren der Mathematischen Wissenschaften, \textbf{266}. Springer-Verlag, Berlin, 1984. 167 pp

\bibitem{Sch_book}
Schwarz, G., \emph{Hodge decomposition---a method for solving boundary value problems},
Lecture Notes in Mathematics, 1607. Springer-Verlag, Berlin, 1995,155 pp.

\bibitem{SylUhl87}
Sylvester, J., Uhlmann, G.,  \emph{A global uniqueness theorem for an inverse boundary value problem},
Ann. of Math. (2) \textbf{125} (1987), no. 1, 153--169.

\bibitem{Tay_bookI}
Taylor, M.,  \emph{Partial differential equations. I. Basic theory}, Applied Mathematical Sciences, 115. Springer-Verlag, New York, 1996. 563 pp.

\bibitem{Tay_bookII}
Taylor, M.,  \emph{Partial differential equations. II. Qualitative studies of linear equations}, Applied Mathematical Sciences, 116. Springer-Verlag, New York, 1996. 528 pp. 




\end{thebibliography}
\end{document}